\documentclass[12pt,twoside,a4paper]{article}
\usepackage{bm}
\usepackage{braket}
\usepackage{graphicx}
\usepackage{graphics}
\usepackage{theorem}
\usepackage{amsmath}
\usepackage{amssymb,mathrsfs}
\usepackage{latexsym}
\usepackage{cases}
\usepackage{bbm}   	%\mathbbm{12}

\newcommand{\one}{\mathbbm{1}}
\usepackage{color}	
\usepackage{algorithm}
\usepackage{algpseudocode}
\usepackage{caption}
\newfloat{procedure}{H}{loa}
\floatname{procedure}{Procedure} %It may be done better
\theorembodyfont{\slshape}
\newtheorem{THM}{Theorem}[section]
\newtheorem{LEM}[THM]{Lemma}%[section]
\newtheorem{PROP}[THM]{Proposition}%[section]
\newtheorem{COR}[THM]{Corollary}%[section]
\newtheorem{ASSU}[THM]{Assumption}

\theorembodyfont{\rmfamily}
\newtheorem{REM}[THM]{Remark}
\newtheorem{rem}[THM]{Remark}
\newtheorem{DEF}[THM]{Definition}

\newenvironment{proof}{\noindent{\it Proof.~~}}{\qed}
\numberwithin{equation}{section}  % If you number theorems, etc. within sections,
\newcommand{\vc}{\bm}
\makeatletter
\DeclareRobustCommand\widecheck[1]{{\mathpalette\@widecheck{#1}}}
\def\@widecheck#1#2{%
    \setbox\z@\hbox{\m@th$#1#2$}%
    \setbox\tw@\hbox{\m@th$#1%
       \widehat{%
          \vrule\@width\z@\@height\ht\z@
          \vrule\@height\z@\@width\wd\z@}$}%
    \dp\tw@-\ht\z@
    \@tempdima\ht\z@ \advance\@tempdima2\ht\tw@ \divide\@tempdima\thr@@
    \setbox\tw@\hbox{%
       \raise\@tempdima\hbox{\scalebox{1}[-1]{\lower\@tempdima\box
\tw@}}}%
    {\ooalign{\box\tw@ \cr \box\z@}}}
\makeatother
\makeatletter
\newif\if@borderstar
\def\bordermatrix{\@ifnextchar*{%
 \@borderstartrue\@bordermatrix@i}{\@borderstarfalse\@bordermatrix@i*}%
}
\def\@bordermatrix@i*{\@ifnextchar[{\@bordermatrix@ii}{\@bordermatrix@ii[()]}}
\def\@bordermatrix@ii[#1]#2{%
\begingroup
 \m@th\@tempdima8.75\p@\setbox\z@\vbox{%
 \def\cr{\crcr\noalign{\kern 2\p@\global\let\cr\endline }}%
 \ialign {$##$\hfil\kern 2\p@\kern\@tempdima & \thinspace %
  \hfil $##$\hfil && \quad\hfil $##$\hfil\crcr\omit\strut %
  \hfil\crcr\noalign{\kern -\baselineskip}#2\crcr\omit %
  \strut\cr}}%
 \setbox\tw@\vbox{\unvcopy\z@\global\setbox\@ne\lastbox}%
 \setbox\tw@\hbox{\unhbox\@ne\unskip\global\setbox\@ne\lastbox}%
 \setbox\tw@\hbox{%
  $\kern\wd\@ne\kern -\@tempdima\left\@firstoftwo#1%
  \if@borderstar\kern 2pt\else\kern -\wd\@ne\fi%
 \global\setbox\@ne\vbox{\box\@ne\if@borderstar\else\kern 2\p@\fi}%
 \vcenter{\if@borderstar\else\kern -\ht\@ne\fi%
  \unvbox\z@\kern -\if@borderstar2\fi\baselineskip}%
 \if@borderstar\kern-2\@tempdima\kern2\p@\else\,\fi\right\@secondoftwo#1 $%
 }\null \;\vbox{\kern\ht\@ne\box\tw@}%
\endgroup
}
\makeatother
\newcommand{\ol}{\overline}
\newcommand{\ool}[1]{\overline{\overline{\bm{#1}}}}
\newcommand{\ul}{\underline}
\newcommand{\wt}{\widetilde}
\newcommand{\wh}{\widehat}
\newcommand{\bv}{\breve}

\newcommand{\down}[2]{\smash{\lower#1\hbox{#2}}}
\newcommand{\up}[2]{\smash{\lower-#1\hbox{#2}}}

\newcommand{\dm}{\displaystyle}
\newcommand{\qed}{\hspace*{\fill}$\Box$}
\newcommand{\vmin}{\wedge}

\newcommand{\vcm}{\ol{\vc{m}}}
\newcommand{\EE}{\mathbb{E}}
\newcommand{\PP}{\mathbb{P}}
\newcommand{\calH}{\mathcal{H}}

\newcommand{\calL}{\mathcal{L}}
\newcommand{\calO}{\mathcal{O}}

\newcommand{\calS}{\mathcal{S}}

\newcommand{\bcalO}{\bm{\mathcal{O}}}
\newcommand{\bbA}{\mathbb{A}}
\newcommand{\bbB}{\mathbb{B}}
\newcommand{\bbC}{\mathbb{C}}

\newcommand{\bbL}{\mathbb{L}}
\newcommand{\bbM}{\mathbb{M}}
\newcommand{\bbN}{\mathbb{N}}
\newcommand{\bbR}{\mathbb{R}}
\newcommand{\bbS}{\mathbb{S}}

\newcommand{\bbZ}{\mathbb{Z}}
\newcommand{\rd}{{\rm d}}

\newcommand{\varep}{\varepsilon}
\makeatletter
\def\eqnarray{\stepcounter{equation}\let\@currentlabel=\theequation
\global\@eqnswtrue
\global\@eqcnt\z@\tabskip\@centering\let\\=\@eqncr
$$\halign to \displaywidth\bgroup\@eqnsel\hskip\@centering
  $\displaystyle\tabskip\z@{##}$&\global\@eqcnt\@ne 
  \hfil$\;{##}\;$\hfil
  &\global\@eqcnt\tw@ $\displaystyle\tabskip\z@{##}$\hfil 
   \tabskip\@centering&\llap{##}\tabskip\z@\cr}
\makeatother
\makeatother

\begin{document}\thispagestyle{empty} 

\hfill
%{\small Last update date: \today}
%Submitted to STOCHASTIC MODELS, \today.

\vspace{-10mm}

{\large{\bf
\begin{center}
A subgeometric convergence formula for finite-level M/G/1-type Markov chains: \\
via a block-decomposition-friendly solution for the Poisson equation of deviation matrix%
\end{center}
}
}

%%%%%%%%%%%%%%%%%%%%%%%%%%%%%%%%%%%%%%%%%%%%%%%%%%%%%%%%%%%%%%%%%
\begin{center}
{
Hiroyuki Masuyama\footnote[1]{E-mail: masuyama@tmu.ac.jp}
}

{\small \it
Graduate School of Management, Tokyo Metropolitan University\\
Tokyo 192--0397, Japan
}

\bigskip
%%%%%%%%%%%%%%%%%%%%%%%%%%%%%%%%%%%%%%%%%%%%%%%%%%%%%%%%%%%%%%%%%
{
Yosuke Katsumata\footnote[2]{E-mail: katsumata@sys.i.kyoto-u.ac.jp}
}

{\small \it
Department of Systems
Science, Graduate School of Informatics, Kyoto University\\
Kyoto 606--8501, Japan
}% \small ends

%%%%%%%%%%%%%%%%%%%%%%%%%%%%%%%%%%%%%%%%%%%%%%%%%%%%%%%%%%%%%%%%%
\bigskip
{
Tatsuaki Kimura%
\footnote[3]{E-mail: kimura@comm.eng.osaka-u.ac.jp}
}

{\small\it
Department of Information and Communications Technology, Graduate School
of Engineering, Osaka University, Suita 565--0871, Japan
}%
%%%%%%%%%%%%%%%%%%%%%%%%%%%%%%%%%%%%%%%%%%%%%%%%%%%%%%%%%%%%%%%%%

\bigskip
\medskip

{\small
\textbf{Abstract}

\medskip

\begin{tabular}{p{0.85\textwidth}}
This paper studies the subgeometric convergence of the stationary distribution in taking the infinite-level limit of a finite-level M/G/1-type Markov chain, that is, in letting the upper boundary level go to infinity. This study is performed through the {\it fundamental deviation matrix}, which is a block-decomposition-friendly solution for the Poisson equation of the deviation matrix. The fundamental deviation matrix yields a difference formula for the respective stationary distributions of the finite-level chain and the corresponding infinite-level chain. The difference formula plays a crucial role in deriving the main result of this paper: a subgeometric convergence formula for the infinite-level limit of the stationary distribution of the finite-level chain. 
\end{tabular}
}
% \samllsize ends
\end{center}

\begin{center}
\begin{tabular}{p{0.90\textwidth}}
{\small

\medskip

}%\samllsize ends
\end{tabular}

\end{center}

%%%%%%%%%%%%%%%%%%%%%%%%%%%%%%%%%%%%%%%%%%%%%%%%%%%%%%%%%%%%%%%%%%%%%
\section{Introduction}\label{introduction}

Finite-level M/G/1-type Markov chains are used for the stationary analysis of finite semi-Markovian queues (see, e.g., \cite{Baio94a,Baio93b,Cuyt03,Herr01}).
However, it is not easy to derive a simple and analytical expression of the stationary distribution in finite-level M/G/1-type Markov chains, except for a few special cases such that the matrix generating functions of level increments are rational \cite{Akar98, KimBara16}.  

There are related studies on the asymptotics of finite-level quasi-birth-and-death processes (QBDs) and finite-level GI/M/1-type Markov chains in taking the infinite-level limit, that is, in letting the upper boundary level go to infinity. As is well known, the QBD is a special case of the finite-level GI/M/1-type Markov chain as well as the finite-level M/G/1-type Markov chain. Miyazawa~et~al.~\cite{Miya07} presented an asymptotic formula for the stationary probability of a finite-level QBD being in the upper boundary level. The asymptotic formula is used to investigate the asymptotic behavior of the loss probability in a MAP/MSP/$c$/$K+c$ queue. J.~Kim and B.~Kim~\cite{KimJeon07} extended the asymptotic formula in \cite{Miya07} to the finite-level GI/M/1-type Markov chain. 

Some researchers have studied the infinite-level limit of finite-level M/G/1-type Markov chains, focusing on the asymptotics of the loss probability in finite M/G/1-type queues, such as queues with Markovian arrival processes.
Ishizaki and Takine~\cite{Ishi99a} established a direct relationship between the respective stationary distributions of a special finite-level M/G/1-type chain and its infinite-level-limit chain (the corresponding infinite-level M/G/1-type chain). Using the direct relationship, the authors obtained the loss probability in a finite M/G/1-type queue with geometrically distributed off-periods. Baiocchi~\cite{Baio94a} derived a geometric asymptotic formula for the loss probability in a MAP/G/1/$K$ queue through the asymptotic analysis of a finite-level M/G/1-type Markov chain with light-tailed level increments. Liu and Zhao \cite{LiuYuan13} presented power-law asymptotic formulas for the loss probability in an M/G/1/$N$ queue with vacations, where the embedded queue length process is a special finite-level M/G/1-type Markov chain with a single background state.

We can find other previous studies \cite{Masu15-AAP,Masu16-SIMAX,Masu17-LAA,Masu17-JORSJ} with related but different perspectives. They are concerned with the upper bound for the error of the  last-column-block-augmented (LCBA) truncation approximation of the stationary distribution in block-structured Markov chains including infinite-level M/G/1-type Markov chains.  Note that the finite-level M/G/1-type Markov chain can be considered the LCBA truncation approximation to an infinite-level M/G/1-type Markov chain. Hence, using the results in those previous studies, we can obtain upper bounds for the difference of the respective stationary distributions of the finite-level and (corresponding) infinite-level M/G/1-type Markov chains.

As we see above, there are no previous studies on the exact convergence speed of the stationary distribution of a finite-level M/G/1-type Markov chain in taking its infinite-level limit. More specifically, no previous studies have presented a convergence formula for the stationary distribution of the finite-level M/G/1-type Markov chain as its upper boundary level goes to infinity.

The main contribution of this paper is to present a {\it subgeometric} convergence formula for the stationary distribution of the finite-level M/G/1-type Markov chain. Subgeometric convergence is much slower than geometric (exponential) convergence, and its typical example is polynomial convergence.

The key to our analysis is a {\it block-decomposition-friendly} solution for the Poisson equation of the {\it deviation matrix} (see, e.g., \cite{Cool02}). We refer to this block-decomposition-friendly solution as the {\it fundamental deviation matrix}. The fundamental deviation matrix yields a difference formula for the respective stationary distributions of the finite-level M/G/1-type Markov chain and its infinite-level-limit chain. Moreover, using the difference formula, we derive a subgeometric convergence formula for the stationary distribution of the finite-level M/G/1-type Markov chain as its upper boundary level goes to infinity. The subgeometric convergence formula requires a certain condition (Assumption~\ref{assumpt-tail}) for the subexponentiality of the integrated tail distribution of nonnegative level increments in steady state of the infinite-level-limit chain.

The rest of this paper consists of five sections. Section~\ref{sec-MG1} describes infinite- and finite-level M/G/1-type Markov chains. Section~\ref{sec-second-order-moment} discusses the second moment condition on level increments of the infinite-level M/G/1-type Markov chain, which leads to the finiteness of the mean of the stationary distribution and to the finiteness of the mean first passage time to level zero. Section~\ref{sec:deviation-matrix} introduces the fundamental deviation matrix of the infinite-level M/G/1-type Markov chain. Based on those results, Section~\ref{sec-convergence} proves the main theorem of this paper, after discussing the uniform convergence of the stationary distribution of the finite-level M/G/1-type Markov chain in taking its infinite-level limit. Finally, Section~\ref{sec:remarks} contains concluding remarks.

%%%%%%%%%%%%%%%%%%%%%%%%%%%%%%%%%%%%%%%%%%%%%%%%%%%%%%%%%%%%%%%%%%%%%%%%%%%%%%%

\section{Model description}\label{sec-MG1}

This section consists of three subsections. Section~\ref{subsec:basic-definition} provides basic definitions and notation. Sections~\ref{subsec:infinite-level-chain} and \ref{subsec:finite-level-chain} describe infinite- and finite-level M/G/1-type Markov chains, respectively.

\subsection{Basic definitions and notation}\label{subsec:basic-definition}

We begin with introducing symbols and notation for numbers. Let
\begin{align*}
\bbZ   &=\{0, \pm1, \pm2,\ldots\}, &
\bbZ_+ &= \{0,1,2,\dots\}, &
\bbN &= \{1,2,3,\dots\},
\\
\bbZ_{\ge k}  &= \{n \in \bbZ: n \ge k\},  & k &\in \bbZ,
\\
\bbZ_{[k,\ell]} &= \{n \in \bbZ: k \le n \le \ell\},  & k,\ell &\in \bbZ,~ k \le \ell,
\end{align*}
and let $\bbM_0$ and $\bbM_1$ denote
\[
\bbM_0 = \bbZ_{[1,M_0]} = \{1, 2, \ldots,M_0\},
\qquad
\bbM_1 = \bbZ_{[1,M_1]} = \{1, 2, \ldots, M_1\},
\]
respectively, where $M_0, M_1 \in \bbN$. For $x,y \in \bbR:=(-\infty,\infty)$, let $x \vmin y = \min(x,y)$ and let $\delta_{k,\ell}$, $k,\ell \in \bbR$, denote the Kronecker delta, that is, $\delta_{k,k} = 1$ and $\delta_{k,\ell} = 0$ for $k \neq \ell$. Furthermore, let $\one(\,\cdot\,)$ denote the indicator function that takes the value of one if the statement in the parentheses is true; otherwise takes the value of zero. 

Next, we describe our notation for vectors and matrices.
All matrices are denoted by bold large letters, and especially, $\vc{O}$ and $\vc{I}$ denote the zero matrix and the identity matrix, respectively, with appropriate sizes (i.e., with appropriate numbers of rows and columns). Basically, all row vectors are denoted by Greek small letters in bold, except for $\vc{g}$ (this exception follows from the convention of the matrix analytic methods pioneered by Neuts \cite{Neut89}); and all column vectors are denoted by English small letters in bold. In particular, $\vc{e}$ denotes the column vector of $1$'s with appropriate sizes. 

Additionally, we introduce more definitions for vectors and matrices. For any matrix (or vector), the absolute value operator $|\,\cdot\,|$ works on it elementwise, and $(\,\cdot\,)_{i,j}$ (resp.\ $(\,\cdot\,)_{i}$) denotes the $(i,j)$-th (resp.~$i$-th) element of the matrix (resp.\ vector) in the parentheses. We then denote by $\|\,\cdot\,\|$, the total-variation norm for vectors. Thus, $\|\vc{s} \| = \sum_{j} |(\vc{s})_j|$ for any vector $\vc{s}$.
Furthermore, for any matrix function $\vc{Z}(\,\cdot\,)$ and scalar function $f(\,\cdot\,)$ on $(-\infty,\infty)$, we use the notations $\vc{Z}(x) = \ol{\bm{\calO}}(f(x)) $ and $\vc{Z}(x) = \ul{\bm{\calO}}(f(x))$:
\begin{align*}
\vc{Z}(x) = \ol{\bm{\calO}}(f(x)) 
&\Longleftrightarrow
\limsup_{x\to\infty} { \sup_i \sum_{j}|(\vc{Z}(x))_{i,j}| \over f(x)} 
< \infty, 
\\
\vc{Z}(x) = \ul{\bm{\calO}}(f(x))
&\Longleftrightarrow
\lim_{x\to\infty}{ \sup_i \sum_{j}|(\vc{Z}(x))_{i,j}| \over f(x)} = 0.
\end{align*}
The notations $\ol{\bm{\calO}}(\,\cdot\,)$ and $\ul{\bm{\calO}}(\,\cdot\,)$ are applied to vector functions; whereas the notations are replaced with $O(\,\cdot\,)$ and $o(\,\cdot\,)$ (according to standard notation) if they are applied to scalar functions. Finally, for any nonnegative matrix $\vc{S} \ge \vc {O}$ (including nonnegative vectors), we write $\vc{S} < \infty$ if every element of $\vc{S}$ is finite.

\subsection{The infinite-level M/G/1-type Markov chain}\label{subsec:infinite-level-chain}

This subsection provides the definition of the infinite-level M/G/1-type Markov chain, and introduces the basic assumption of the present paper.

We describe the infinite-level M/G/1-type Markov chain. Let $\{(X_n, J_n); n \in \bbZ_+\}$ denote a discrete-time Markov chain on state space $\bbS: = \bigcup_{k=0}^{\infty} \bbL_k$, where $\bbL_k = \{k\} \times \bbM_{k \vmin 1}$ for $k \in \bbZ_+$. The subset $\bbL_k$ of state space $\bbS$ is referred to as {\it level $k$}. Let $\vc{P}$ denote the transition probability matrix of the Markov chain $\{(X_n, J_n)\}$, and assume that $\vc{P}$ is a stochastic matrix such that
\begin{equation}
\vc{P}
=\bordermatrix{
&  \bbL_0         
& \bbL_1   
& \bbL_2 
& \bbL_3   
& \cdots 
\cr
	\bbL_0
& \bm{B}(0)        
& \bm{B}(1) 
& \bm{B}(2) 
& \bm{B}(3) 
& \cdots 
\cr
	\bbL_1
& \bm{B}(-1)            
& \bm{A}(0) 
& \bm{A}(1) 
& \bm{A}(2)  
& \cdots 
\cr
	\bbL_2
& \bm{O} 
& \bm{A}(-1) 
& \bm{A}(0) 
& \bm{A}(1)  
& \cdots 
\cr
	\bbL_3
& \bm{O} 
& \bm{O} 
& \bm{A}(-1) 
& \bm{A}(0)  
& \cdots 
\cr
~\vdots  
& \vdots
& \vdots
& \vdots   
& \vdots     
& \ddots
}.
\label{defn-P}
\end{equation}
Therefore, the component block matrices $\vc{A}(k)$ and $\vc{B}(k)$ satisfy the following:
\begin{align}
\sum_{k=-1}^{\infty}\vc{A}(k)\vc{e} &= \vc{e},
\qquad
\sum_{k=0}^{\infty}\vc{B}(k)\vc{e} = \vc{e},
\label{cond-finite-MG1-01}
\\
\vc{B}(-1)\vc{e} &= \vc{A}(-1)\vc{e}.
\label{cond-finite-MG1-03}
\end{align}
%.
The Markov chain $\{(X_n, J_n)\}$ is referred to as an {\it M/G/1-type Markov chain} (see \cite{Neut89}). We call this Markov chain the {\it infinite-level M/G/1-type Markov chain} or {\it infinite-level chain} for short, in order to distinguish it from its finite-level version (described later).

We introduce our basic assumption. To this end, let 
\begin{eqnarray}
\vc{A} &=& \sum_{k=-1}^{\infty}\vc{A}(k),
\qquad
\vcm_{A} = \sum_{k=-1}^{\infty} k \vc{A}(k)\vc{e},
\label{defn-beta_A}
\end{eqnarray}
where $\vc{A}$ is a stochastic matrix due to (\ref{cond-finite-MG1-01}). The following is then the basic assumption.
\begin{ASSU} \label{assumpt-ergodic-MG1}
(i) The stochastic matrices $\vc{A}$ and $\vc{P}$ (given in (\ref{defn-P})) are irreducible; 
(ii) $\vcm_{B}:=\sum_{k=1}^{\infty} k\vc{B}(k)\vc{e} < \infty$; and 
(iii) $\sigma :=\vc{\varpi}\vcm_{A} < 0$, 
where $\vc{\varpi}$ denotes the unique stationary distribution vector of~$\vc{A}$.
\end{ASSU}

Assumption~\ref{assumpt-ergodic-MG1} ensures (see, e.g., \cite[Chapter XI, Proposition~3.1]{Asmu03}) that the infinite-level chain $\{(X_n, J_n)\}$ is irreducible and positive recurrent and therefore this chain has the unique stationary distribution vector, denoted by $\vc{\pi} = (\pi(k,i))_{(k,i)\in\bbS}$. 
For later convenience, $\vc{\pi}$ is partitioned level-wise: $\vc{\pi} = (\vc{\pi}(0),\vc{\pi}(1),\dots)$, where $\vc{\pi}(k) = (\pi(k,i))_{i \in \bbM_{k \vmin 1}}$ for $k \in \bbZ_+$.

To describe the stationary distribution $\vc{\pi}$, we introduce the $G$- and $R$-matrices of the infinite-level M/G/1-type Markov chain. Let $\vc{G}:=(G_{i,j})_{i,j\in\bbM_1}$ denote an $M_1 \times M_1$ matrix such that
\[
G_{i,j}
=
\PP (J_{\tau_k} = j \mid (X_0, J_0) = (k+1, i) \in \bbL_{\ge 2}),
\]
where $\tau_k = \inf\{n \in \bbN: X_n = k\}$ and $\bbL_{\ge k}= \bigcup_{\ell=k}^{\infty} \bbL_{\ell}$ for $k \in \bbZ_+$. 
Assumption~\ref{assumpt-ergodic-MG1} (i) and (iii) ensures that $\vc{G}$ is a stochastic matrix with a single closed communicating class \cite[Proposition 2.1]{Kimu10} and thus the unique stationary distribution vector, denoted by $\vc{g}$.
Additionally, using the $G$-matrix $\vc{G}$, we define $\vc{R}_0(k)$ and
$\vc{R}(k)$, $k \in \bbN$, as
\begin{alignat}{2}
\vc{R}_0(k)
&=
\sum_{m=k}^{\infty}
\vc{B}(m)\vc{G}^{m-k} (\vc{I}-\vc{\Phi}(0))^{-1},& \qquad k & \in \bbN,
\label{def-R_0(k)}
\\
\vc{R}(k)
&=
\sum_{m=k}^{\infty}
\vc{A}(m)\vc{G}^{m-k} (\vc{I}-\vc{\Phi}(0))^{-1},& \qquad k & \in \bbN,
\label{def-R(k)}
\end{alignat}
respectively, where
\begin{equation}
\vc{\Phi}(0)
= \sum_{m=0}^{\infty} \vc{A}(m)\vc{G}^{m}.
\label{defn-Phi(0)}
\end{equation}
We then have
\begin{equation}
\vc{\pi}(k) = \vc{\pi}(0)\vc{R}_0(k)
+ \sum_{\ell=1}^{k-1}\vc{\pi}(\ell)\vc{R}(k-\ell),
\qquad k \in \bbN,
\label{eqn-pi(k)}
\end{equation}
which is referred to as {\it Ramaswami's recursion} \cite{Rama88}.

\subsection{The finite-level M/G/1-type Markov chain}\label{subsec:finite-level-chain}

In this subsection, we define the finite-level M/G/1-type Markov chain and its stationary distribution vector. We also introduce our convention for handling calculations over finite and infinite matrices (including vectors). The convention facilitates the discussion in the subsequent sections. In addition, symbol $N$ is assumed to an arbitrary positive integer throughout the paper, unless otherwise stated.

We provide the definition of the finite-level M/G/1-type Markov chain. For any $N \in \bbN$, let $\{(X_n^{(N)}, J_n^{(N)}); n\in\bbZ_+\}$ denote a discrete-time Markov chain on state space $\bbL_{\le N}: = \bigcup_{\ell=0}^N \bbL_{\ell}$ with transition probability matrix $\vc{P}^{(N)}$ given by
{\let\quad\thinspace
\begin{equation}
\vc{P}^{(N)} =
\bordermatrix{
& \bbL_0
& \bbL_1
& \bbL_2
& \cdots
& \bbL_{N-2}
& \bbL_{N-1}
& \bbL_N
\cr
\bbL_0
& \vc{B}(0) 
& \vc{B}(1) 
& \vc{B}(2) 
& \cdots 
& \vc{B}(N-2) 
& \vc{B}(N-1) 
& \ol{\vc{B}}(N-1)  
\cr
\bbL_1
& \vc{B}(-1) 
& \vc{A}(0) 
& \vc{A}(1) 
& \cdots 
& \vc{A}(N-3) 
& \vc{A}(N-2) 
& \ol{\vc{A}}(N-2) 	 
\cr
\bbL_2
& \vc{O} 
& \vc{A}(-1) 
& \vc{A}(0) 
& \cdots	
& \vc{A}(N-4) 
& \vc{A}(N-3) 
& \ol{\vc{A}}(N-3)  	
\cr
~\vdots
& \vdots 
& \vdots 
& \vdots 
& \ddots 
& \vdots 
& \vdots 
& \vdots
\cr
\bbL_{N-2}
& \vc{O} 
& \vc{O} 
& \vc{O} 
& \cdots 
& \vc{A}(0) 
& \vc{A}(1) 
& \ol{\vc{A}}(1) 
\cr
\bbL_{N-1}
& \vc{O} 
& \vc{O} 
& \vc{O} 
& \cdots 
& \vc{A}(-1) 
& \vc{A}(0) 
& \ol{\vc{A}}(0) 
\cr
\bbL_N
& \vc{O} 
& \vc{O} 
& \vc{O} 
& \cdots 
& \vc{O} 
& \vc{A}(-1) 
& \ol{\vc{A}}(-1) 
},
\label{defn:P^{(N)}}
\end{equation}
}
where
\begin{subequations}\label{eqn:ol{A}(k)-ol{B}(k)}
\begin{alignat}{2}
\ol{\vc{A}}(k)
&= \sum_{\ell=k+1}^{\infty} \vc{A}(\ell), 
&\qquad
k &\in \bbZ_{\ge -2}, 
\label{defn-ol{A}(k)}
\\
\ol{\vc{B}}(k)
&= \sum_{\ell=k+1}^{\infty} \vc{B}(\ell), 
&\qquad
k &\in \bbZ_+.
\label{defn-ol{B}(k)}
\end{alignat}
\end{subequations}
The Markov chain $\{(X_n^{(N)}, J_n^{(N)})\}$ is referred to as the {\it finite-level M/G/1-type Markov chain} or {\it finite-level chain} for short.  Without loss of generality, we assume that the finite-level chain $\{(X_n^{(N)}, J_n^{(N)})\}$ is defined on the same probability space as that of the infinite-level chain $\{(X_n,J_n)\}$.
\begin{REM}
The stochastic matrix $\vc{P}^{(N)}$ can be considered the {\it last-column-block-augmented (LCBA) truncation} of $\vc{P}$ (see \cite{Masu15-AAP,Masu16-SIMAX,Masu17-LAA,Masu17-JORSJ}).
\end{REM}

Proposition~\ref{prop:P^{(N)}} below is concerned with the recurrence of $\{(X_n^{(N)}, J_n^{(N)})\}$.
\begin{PROP}\label{prop:P^{(N)}}
Under Assumption~\ref{assumpt-ergodic-MG1}, the following hold for level zero $\bbL_0$ of the finite-level chain $\{(X_n^{(N)}, J_n^{(N)})\}$:
\begin{enumerate}
\item For all sufficiently large $N \in \bbN$, any state $(0,j) \in \bbL_0$ is accessible from any state $(0,i) \in \bbL_0$.
\item For each $N \in \bbN$, level zero $\bbL_0$ is accessible from any state in its state space $\bbL_{\le N}$ with probability one.

\end{enumerate}
\end{PROP}

\begin{proof}
We begin with the proof of the statement (i). Assumption~\ref{assumpt-ergodic-MG1} ensures that the infinite-level chain $\{(X_n, J_n)\}$ is irreducible and thus it reaches any state $(0,j) \in \bbL_0$ from any state $(0,i) \in \bbL_0$. Let $N_{i,j}$ denote the maximum level in an arbitrarily chosen path from state $(0,i)$ to state $(0,j)$. Since $\bbL_0$ is finite, we have $N_0:=\max_{i,j \in \bbM_0} N_{i,j} < \infty$. Furthermore, the similarity in the structures of $\vc{P}^{(N)}$ and $\vc{P}$ implies that the finite-level chain $\{(X_n^{(N)}, J_n^{(N)})\}$ behaves stochastically the same as the infinite-level chain $\{(X_n, J_n)\}$ unless the former reaches its upper boundary level $N$. Therefore, for any $N > N_0$, the finite-level chain $\{(X_n^{(N)}, J_n^{(N)})\}$ reaches any state $(0,j) \in \bbL_0$ from any state $(0,i) \in \bbL_0$.

Next, we prove the statement (ii). To do this, we construct the sample paths of the finite-level chain $\{(X_n^{(N)}, J_n^{(N)})\}$ before the first passage time to level zero by using the sample paths of the infinite-level chain $\{(X_n, J_n)\}$, based on Proposition~\ref{PROP:Appendix-01}:
\begin{subequations}\label{defn:(X_n^{(N)},J_n^{(N)})}
\begin{alignat}{2}
X_n^{(N)} &= 
\left\{
\begin{array}{@{\ }l@{\ }l}
X_0 \vmin N, & n=0,
\\
(X_{n-1}^{(N)} + X_n - X_{n-1}) \vmin N, &
n=1,2,\dots,\tau_{0}^{(N)},~X_{n-1}^{(N)} \in \bbZ_{[0,N-1]},
\\
X_{n-1}^{(N)} + (X_n - X_{n-1}) \vmin 0, &
n=1,2,\dots,\tau_{0}^{(N)},~X_{n-1}^{(N)} =N,
\end{array}
\right.
\label{defn:X_n^{(N)}}
\\
J_n^{(N)} &= J_n, \qquad n=0,1,\dots,\tau_{0}^{(N)}-1,
\label{defn:J_n^{(N)}}
\end{alignat}
\end{subequations}
where $\tau_0^{(N)} = \inf\{n \in \bbN:X_n^{(N)} = 0\}$. The first and second cases of (\ref{defn:X_n^{(N)}}) (together with (\ref{defn:J_n^{(N)}})) ensure that the level increment of $\{(X_n^{(N)}, J_n^{(N)})\}$ is equal to that of $\{(X_n, J_n)\}$ while $\{(X_n^{(N)}, J_n^{(N)})\}$ is below level $N$ and does not reach level zero yet. In addition, the third case of (\ref{defn:X_n^{(N)}}) ensures that $\{(X_n^{(N)}, J_n^{(N)})\}$ stays at level $N$ while $\{(X_n, J_n)\}$ evolves with the transition matrices $\vc{A}(0), \vc{A}(1), \vc{A}(2),\dots$; and ensures that $\{(X_n^{(N)}, J_n^{(N)})\}$ moves down to level $N-1$ when $\{(X_n, J_n)\}$ moves down by one level with $\vc{A}(-1)$. Thus, the behavior of $\{(X_n^{(N)}, J_n^{(N)})\}$ specified by (\ref{defn:(X_n^{(N)},J_n^{(N)})}) follows the transition law (\ref{defn:P^{(N)}}).

The pass construction (\ref{defn:(X_n^{(N)},J_n^{(N)})}) contributes to the proof of the statement (i). Indeed, (\ref{defn:(X_n^{(N)},J_n^{(N)})}) implies that $X_n^{(N)} \le X_n$ and $J_n^{(N)} = J_n$ for all $n = 0,1,\dots,\tau_0^{(N)}-1$. In addition, Assumption~\ref{assumpt-ergodic-MG1} ensures that the infinite-level chain $\{(X_n, J_n)\}$ is irreducible and positive recurrent.
Therefore, the finite-level chain $\{(X_n^{(N)}, J_n^{(N)})\}$ reaches $\bbL_0$ from all the states in $\bbL_{\le N}$ with probability one.
\end{proof}

The finite-level chain $\{(X_n^{(N)}, J_n^{(N)})\}$ has at least one stationary distribution vector, denoted by $\vc{\pi}^{(N)}:=(\pi^{(N)}(k,i))_{(k,i)\in\bbL_{\le N}}$, which is uniquely determined for all sufficiently large $N \in \bbN$ due to Proposition~\ref{prop:P^{(N)}}. By definition,
\begin{equation}
\vc{\pi}^{(N)}\vc{P}^{(N)} = \vc{\pi}^{(N)},
\quad \vc{\pi}^{(N)}\vc{e} = 1,
\quad \vc{\pi}^{(N)} \ge \vc{0}.
\label{defn-pi^{(N)}}
\end{equation}
For later use, $\vc{\pi}^{(N)}$ is partitioned as 
\[
\vc{\pi}^{(N)} = (\vc{\pi}^{(N)}(0),\vc{\pi}^{(N)}(1),\dots,\vc{\pi}^{(N)}(N)),
\]
where $\vc{\pi}^{(N)}(k) = (\pi^{(N)}(k,i))_{i \in \bbM_{k \vmin 1}}$ for $k \in \bbZ_{[0,N]}$. 

Finally, we introduce our convention for performing calculations (such as addition and multiplication) over finite and infinite dimensional matrices (including vectors). As mentioned in the introduction, the main purpose of this paper is to study the convergence of $\{\vc{\pi}^{(N)}; N \in\bbN\}$. Thus, we consider the situation that the probability vectors $\vc{\pi}^{(N)}$, $N \in\bbN$ of different finite dimensions converge to 
a certain probability vector of infinite dimension (which is equal to $\vc{\pi}$, as expected). To facility this study, the following rule is introduced: Finite dimensional matrices (including vectors) are extended (if necessary) to infinite dimensional matrices by appending zeros to them keeping their original elements in the original positions. According to this rule, for example, $\vc{\pi}^{(N)} - \vc{\pi}$ and $\vc{P}^{(N)} - \vc{P}$ are well-defined.

\section{Second-order moment condition on level increments in the infinite-level chain}\label{sec-second-order-moment}

This section introduces a condition on level increments of the infinite-level chain (Assumption~\ref{assum-2nd-moment} below). For reference, we call it the {\it second-order moment condition}. As shown later, the second-order moment condition ensures that the stationary distribution is finite and that the mean first passage time to level zero is finite, which leads to the proof of the convergence of $\{\vc{\pi}^{(N)}\}$ to $\vc{\pi}$.
\begin{ASSU}[Second-order moment condition on level increments] \label{assum-2nd-moment}
\[
\sum_{k=1}^{\infty} k^2 \vc{A}(k) < \infty,
\qquad
\sum_{k=1}^{\infty} k^2 \vc{B}(k) < \infty.
\]
\end{ASSU}

\begin{REM}\label{rem:assum-2nd-moment}
Assumption~\ref{assum-2nd-moment} implies that $\vc{A}(k) = \ul{\bcalO}(k^{-3})$ and $\vc{B}(k) = \ul{\bcalO}(k^{-3})$.
\end{REM}

\subsection{Finiteness of the mean of the stationary distribution}

In this subsection, we first establish a certain Foster-Lyapunov drift condition (called a drift condition for short). Using this drift condition, we show that the second-order moment condition is equivalent to $\sum_{k=1}^{\infty}k\vc{\pi}(k)\vc{e} < \infty$, and also show that the second-order moment condition implies $\sup_{N \in \bbN}\sum_{k=1}^N k \vc{\pi}^{(N)}(k) \vc{e} < \infty$.

We consider the Poisson equation
\begin{equation}
(\vc{I} - \vc{A})\vc{x} 
= -\sigma\vc{e} + \vcm_{A}
\label{Poisson-EQ-A}
\end{equation}
to establish the drift condition under Assumption~\ref{assum-2nd-moment}. 
A solution for this Poisson equation is given by
\begin{equation}
\vc{a} = (\vc{I} - \vc{A} + \vc{e}\vc{\varpi})^{-1}\vcm_{A}
+ c \vc{e},
\label{defn-alpha}
\end{equation}
where $c\in (-\infty,\infty)$ is an arbitrary constant. 
Using $\vc{A}\vc{e}=\vc{e}$, $\sigma=\vc{\varpi}\vcm_{A}$, and $\vc{\varpi}(\vc{I} - \vc{A} + \vc{e}\vc{\varpi})^{-1}=\vc{\varpi}$, we can readily confirm that 
\begin{align}
(\vc{I} - \vc{A})\vc{a} 
= -\sigma\vc{e} + \vcm_{A}.
\label{eqn:a}
\end{align}
For later use, we fix $c > 0$ sufficiently large such that $\vc{a} \ge \vc{0}$, and then define $\vc{v}:=(v(k,i))_{(k,i)\in\bbS}$ and $\vc{f}:=(f(k,i))_{(k,i)\in\bbS}$ as nonnegative column vectors such that
\begin{alignat}{2}
\vc{v}(k):= (v(k,i))_{i\in\bbM_{k\vmin1}}
&= 
\left\{
\begin{array}{ll}
\vc{0}, & k =0,
\\
\dm{1 \over -\sigma} \left(k^2\vc{e} + 2k\vc{a} \right), & k \in\bbN,
\end{array}
\right.
\label{defn-v(k)}
\\
\vc{f}(k):= (f(k,i))_{i\in\bbM_{k\vmin1}}
&=\left\{
\begin{array}{ll}
\vc{e}, & k=0,
\\
(k+1)\vc{e}, & k \in\bbN,
\end{array}
\right.
\label{defn-f(k)}
\end{alignat}
respectively, where the size $M_0$ of $\vc{v}(0)$ and $\vc{f}(0)$ is, in general, different from the size $M_1$ of $\vc{v}(k)$ and $\vc{f}(k)$, $k \in \bbN$. Furthermore, let $\vc{1}_{\bbC}:=(1_{\bbC}(k,i))_{(k,i)\in\bbS}$, $\bbC \subseteq \bbS$, denote a column vector such that
\[
1_{\bbC}(k,i) = 
\left\{
\begin{array}{ll}
1, & (k,i) \in \bbC,
\\
0, & (k,i) \not\in \bbC.
\end{array}
\right.
\]
When $\bbC$ is a single state $(\ell,j) \in \bbS$, that is, $\bbC = \{(\ell,j)\}$, we write $\vc{1}_{(\ell,j)}$ for $\vc{1}_{\{(\ell,j)\}}$.

The following lemma presents our desired drift condition.
\begin{LEM}\label{lem-drift-cond}
If Assumptions~\ref{assumpt-ergodic-MG1} and \ref{assum-2nd-moment} hold, then there exist some $b \in (0,\infty)$ and $K \in \bbN$ such that
\begin{align}
\vc{P}^{(N)}\vc{v}
\le
\vc{P}\vc{v}
\le \vc{v} - \vc{f} + b \vc{1}_{\bbL_{\le K}}\quad \mbox{for all $N \in \bbN$},
\label{drift-cond-Pv}
\end{align}
where $\bbL_{\le k} = \bigcup_{\ell=0}^k \bbL_{\ell}$ for $k \in \bbZ_+$. 
\end{LEM}
\begin{proof}
To prove this lemma, it suffices to show that
\begin{subequations}\label{drift-cond-Pv-02}
\begin{alignat}{2}
\sum_{\ell=0}^{\infty}\vc{P}(k;\ell)\vc{v}(\ell) &< \infty
& &\mbox{for all $k \in \bbZ_+$},
\label{ineqn-sum-P(k;l)v(l)-a}
\\
\sum_{\ell=0}^{\infty}\vc{P}(k;\ell)\vc{v}(\ell)
&\le \vc{v}(k) - \vc{f}(k) &\quad &\mbox{for all sufficiently large $k \in \bbZ_+$},
\label{ineqn-sum-P(k;l)v(l)-b}
\end{alignat}
\end{subequations}
where $\vc{P}(k;\ell)$ $k,\ell \in \bbZ_+$, denotes a submatrix of $\vc{P}$ that contains the transition probabilities from level $k$ to level $\ell$. 
Indeed, (\ref{drift-cond-Pv-02}) implies that there exist some $b \in (0,\infty)$ and $K \in \bbN$ such that
\begin{equation*}
\vc{P}\vc{v}
\le \vc{v} - \vc{f} + b \vc{1}_{\bbL_{\le K}}.
\end{equation*}
Furthermore, $\vc{P}^{(N)}\vc{v} \le \vc{P}\vc{v}$ for $N \in \bbN$. This inequality follows from (\ref{defn-P}) and (\ref{defn:P^{(N)}}) together with that $\vc{v}(1) \le \vc{v}(2) \le \vc{v}(3) \le \cdots$ due to (\ref{defn-v(k)}). 

We prove (\ref{ineqn-sum-P(k;l)v(l)-a}). It follows from (\ref{defn-P}), (\ref{defn-beta_A}), and (\ref{defn-v(k)}) that, for all $k \in \bbZ_{\ge 2}$,
\begin{align}
&
\sum_{\ell=0}^{\infty}\vc{P}(k;\ell)\vc{v}(\ell)
\nonumber
\\
&\quad =
\sum_{\ell=-1}^{\infty} \vc{A}(\ell)\vc{v}(k+\ell)
=
{1 \over -\sigma} 
\left[
\sum_{\ell=-1}^{\infty} (k+\ell)^2 \vc{A}(\ell)\vc{e}
+ 2\sum_{\ell=-1}^{\infty} (k+\ell) \vc{A}(\ell)\vc{a}
\right]
\nonumber
\\
&\quad =
{1 \over -\sigma} 
\left[ k^2 \vc{e}
+ 2 k \left( 
\vcm_{A} + \vc{A}\vc{a}
\right)
\right]
+   {1 \over -\sigma} 
\left[
\sum_{\ell=-1}^{\infty} \ell^2 \vc{A}(\ell)\vc{e}
+  2\sum_{\ell=-1}^{\infty} \ell \vc{A}(\ell)\vc{a}
\right]
\nonumber
\\
&\quad =
{1 \over -\sigma} 
\left[
k^2 \vc{e} + 2k ( \vc{a} + \sigma\vc{e} )
\right]
+  
{1 \over -\sigma} 
\left[
\sum_{\ell=-1}^{\infty} \ell^2 \vc{A}(\ell)\vc{e}
+  2\sum_{\ell=-1}^{\infty} \ell \vc{A}(\ell)\vc{a}
\right],
\label{eqn-180628-01}
\end{align}
where the last equality holds due to (\ref{eqn:a}).
It also follows from (\ref{eqn-180628-01}) and Assumption~\ref{assum-2nd-moment} that $\sum_{\ell=0}^{\infty}\vc{P}(k;\ell)\vc{v}(\ell)$ is finite for each $k \in \bbZ_{\ge 2}$. Similarly, we can confirm that $\sum_{\ell=0}^{\infty}\vc{P}(0;\ell)\vc{v}(\ell)$ and $\sum_{\ell=0}^{\infty}\vc{P}(1;\ell)\vc{v}(\ell)$ are finite. 

Next, we prove (\ref{ineqn-sum-P(k;l)v(l)-b}). Using (\ref{defn-v(k)}) and (\ref{defn-f(k)}), we rewrite (\ref{eqn-180628-01}) as
\begin{align}
&\sum_{\ell=0}^{\infty}\vc{P}(k;\ell)\vc{v}(\ell)
\nonumber
\\
&~
= \vc{v}(k) - 2k\vc{e}
+ {1 \over -\sigma} 
\left[
\sum_{\ell=-1}^{\infty} \ell^2 \vc{A}(\ell)\vc{e}
+  2\sum_{\ell=-1}^{\infty} \ell \vc{A}(\ell)\vc{a}
\right]
\nonumber
\\
&~= \vc{v}(k) - \vc{f}(k) - (k-1)\vc{e}
+
{1 \over -\sigma}
\left[
 \sum_{\ell=-1}^{\infty} \ell^2 \vc{A}(\ell)\vc{e}
+  2\sum_{\ell=-1}^{\infty} \ell \vc{A}(\ell)\vc{a}
\right],~~ k \in \bbZ_{\ge 2}.
\label{eqn-180628-02}
\end{align}
Clearly, there exists some $K \in \bbN$ such that
\[
- (k-1)\vc{e}
+
{1 \over -\sigma}
\left[
 \sum_{\ell=-1}^{\infty} \ell^2 \vc{A}(\ell)\vc{e}
+  2\sum_{\ell=-1}^{\infty} \ell \vc{A}(\ell)\vc{a}
\right]
\le \vc{0}\quad \mbox{for all $k \in \bbZ_{\ge K+1}$}.
\]
Combining this and (\ref{eqn-180628-02}) results in (\ref{ineqn-sum-P(k;l)v(l)-b}). The proof is completed. %\qed
\end{proof}

Lemma~\ref{lem-drift-cond} leads to the results on the finite means of the stationary distribution vectors $\vc{\pi}$ and $\vc{\pi}^{(N)}$ of the infinite-level and finite-level chains.
\begin{THM}\label{thm-mean}
Under Assumption~\ref{assumpt-ergodic-MG1}, the following are true:
\begin{enumerate}
\item Assumption~\ref{assum-2nd-moment} 
holds if and only if
\begin{eqnarray}
\sum_{k=1}^{\infty} k \vc{\pi}(k) \vc{e} &<& \infty.
\label{sum-k*pi(k)e<b}
\end{eqnarray}
\item If Assumption~\ref{assum-2nd-moment} holds, then
\begin{eqnarray}
\sup_{N \in \bbN}\sum_{k=1}^N k \vc{\pi}^{(N)}(k) \vc{e} &<& \infty.
\label{sum-k*pi^{(N)}(k)e<b}
\end{eqnarray}
\end{enumerate}
\end{THM}

\begin{proof}
The ``only if" part of the statement (i), that is, ``Assumption~\ref{assum-2nd-moment}  $\Longrightarrow$ (\ref{sum-k*pi(k)e<b})", can be proved in a similar way to the proof of the statement (ii). Thus in what follows, we prove the statement (ii) and the ``if" part of the statement (i). 

We prove the statement (ii). Suppose that Assumption~\ref{assum-2nd-moment} holds (in addition to Assumption~\ref{assumpt-ergodic-MG1}) and thus Lemma~\ref{lem-drift-cond} holds. Pre-multiplying the left- and right-hand sides of (\ref{drift-cond-Pv}) by $\vc{\pi}^{(N)}$ and using (\ref{defn-pi^{(N)}}), we obtain
\begin{align*}
\vc{\pi}^{(N)}\vc{v}
&= \vc{\pi}^{(N)} \vc{P}^{(N)}\vc{v}
\le \vc{\pi}^{(N)} \vc{v} - \vc{\pi}^{(N)}\vc{f} + b,
\end{align*}
and thus $\vc{\pi}^{(N)}\vc{f} \le b$ for all $N \in \bbN$.  Combining this inequality and (\ref{defn-f(k)}) yields (\ref{sum-k*pi^{(N)}(k)e<b}). The statement (ii) has been proved.

We prove the ``if" part of the statement (i). To this end, suppose that (\ref{sum-k*pi(k)e<b}) holds. It then follows from (\ref{eqn-pi(k)}) that
\begin{align*}
\infty 
> \sum_{k=1}^{\infty} k\vc{\pi}(k)\vc{e} 
&=\vc{\pi}(0) \sum_{k=1}^{\infty} k\vc{R}_0(k)\vc{e} 
+ 
\sum_{\ell=1}^{\infty} \vc{\pi}(\ell)  
\sum_{k=\ell+1}^{\infty} k \vc{R}(k-\ell)\vc{e}
\nonumber
\\
&=\vc{\pi}(0) \sum_{k=1}^{\infty} k\vc{R}_0(k)\vc{e} 
+ 
\sum_{\ell=1}^{\infty} \vc{\pi}(\ell)  
\sum_{k=\ell+1}^{\infty} (k-\ell) \vc{R}(k-\ell)\vc{e}
\nonumber
\\
&{} \quad + 
\sum_{\ell=1}^{\infty} \ell \vc{\pi}(\ell)  
\sum_{k=\ell+1}^{\infty} \vc{R}(k-\ell)\vc{e},
%\label{eqn-sum-k*pi(k)}
\end{align*}
which yields $\sum_{k=1}^{\infty} k\vc{R}(k)\vc{e} < \infty$ and $\sum_{k=1}^{\infty} k\vc{R}_0(k)\vc{e}< \infty$. It also follows from (\ref{def-R(k)}), $\vc{G}\vc{e}=\vc{e}$, and $(\vc{I} - \vc{\Phi}(0))^{-1}\vc{e} \ge \vc{e}$ that
\begin{eqnarray*}
\infty > \sum_{k=1}^{\infty} k\vc{R}(k)\vc{e} 
&=&
\sum_{k=1}^{\infty} k
\sum_{\ell=k}^{\infty}
\vc{A}(\ell)\vc{G}^{\ell-k} (\vc{I}-\vc{\Phi}(0))^{-1}\vc{e}
\nonumber
\\
&\ge&
\sum_{k=1}^{\infty} k
\sum_{\ell=k}^{\infty}
\vc{A}(\ell)\vc{G}^{\ell-k}\vc{e}
= \sum_{k=1}^{\infty} k
\sum_{\ell=k}^{\infty}
\vc{A}(\ell)\vc{e}
\nonumber
\\
&=&
{1 \over 2}\sum_{\ell=1}^{\infty} \ell(\ell+1)
\vc{A}(\ell)\vc{e},
\end{eqnarray*}
which shows that $\sum_{\ell=1}^{\infty} \ell^2\vc{A}(\ell)\vc{e} < \infty$.
Similarly, combining (\ref{def-R_0(k)}) and $\sum_{k=1}^{\infty} k\vc{R}_0(k)\vc{e} < \infty$ leads to $\sum_{\ell=1}^{\infty} \ell^2\vc{B}(\ell)\vc{e} < \infty$. Consequently, the ``if" part of the statement (i) has been proved. 
\end{proof}

\subsection{Finiteness of the mean first passage time to level zero}\label{subsec:mean_first passage_time}

This subsection provides basic results on the mean first passage time to level zero of the infinite-level chain. First, we introduce some definitions to derive the results. We then show that the mean first passage time to level zero is basically linear with the starting level. We also show that the second-order moment condition (Assumption~\ref{assum-2nd-moment}) is equivalent to that the mean first passage time to level zero in steady state is finite. These results are used in the subsequent sections. Although some of them might be found in the literature, we here provide them for the reader's convenience.

We begin with some definitions associated with the mean first passage time to level zero. Let
$\vc{u}(k):=(u(k,i))_{i\in\bbM_{k\vmin1}}$, $k\in\bbZ_+$, denote a column vector such that
\begin{equation}
u(k,i) = \EE_{(k,i)}[\tau_0] \ge 1,\qquad (k,i) \in \bbS,
\label{defn-u(k,i)}
\end{equation}
where $\tau_0 = \inf\{n \in \bbN: X_n = k\}$ for $k \in \bbZ_+$. 
Let $\wh{\vc{G}}_1(z)$, $z \in [0,1]$, denote an $M_1 \times M_0$ matrix such that
\begin{alignat}{2}
(\wh{\vc{G}}_1(z))_{i,j}
&= \EE_{(1,i)} \!\left [z^{\tau_0}\one(J_{\tau_0} = j) \right],
\qquad i \in \bbM_1,~ j \in \bbM_0,
\label{defn:wh{G}_{1,i,j}(z)}
\end{alignat}
where $\EE_{(k,i)}[\,\cdot\,] = \EE[\,\cdot\mid (X_0, J_0)=(k,i)]$ for $(k,i) \in \bbS$. Furthermore, let $\wh{\vc{G}}(z)$, $z \in [0,1]$, denote an $M_1 \times M_1$ matrix such that
\begin{alignat*}{2}
(\wh{\vc{G}}(z))_{i,j}
&= \EE_{(2,i)}\!\left [z^{\tau_1}\one(J_{\tau_1} = j)\right],
\qquad i,j \in \bbM_1.
%\label{defn:wh{G}(z)}
\end{alignat*}
Note here that, except for level zero, the infinite-level M/G/1-type Markov chain has the level homogeneity of transitions. Therefore, for all $k\in\bbN$,
\begin{alignat}{2}
(\wh{\vc{G}}(z))_{i,j}
&= \EE_{(k+1,i)}\!\left [z^{\tau_k}\one(J_{\tau_k} = j)\right],
\qquad i,j \in \bbM_1,
\label{defn:wh{G}(z)-org}
\end{alignat}
and $\wh{\vc{G}}(1)$ is equal to the $G$-matrix $\vc{G}$.  In addition, it follows from (\ref{defn-u(k,i)}), (\ref{defn:wh{G}_{1,i,j}(z)}), and (\ref{defn:wh{G}(z)-org}) that
\begin{eqnarray}
\vc{u}(k) 
&=& \left. {\rd \over \rd z} [ \wh{\vc{G}}(z) ]^{k-1}\wh{\vc{G}}_1(z) \right|_{z=1} \vc{e},\qquad k \in \bbN.
\label{eqn-u(k)-02}
\end{eqnarray}

The following lemma shows that the mean first passage time to zero from level $k$ has its dominant term linear with $k$.
\begin{LEM}\label{lem-u(k)}
Suppose that Assumption~\ref{assumpt-ergodic-MG1} holds. We then have
\begin{subnumcases}
{\label{eqn-u}}
\vc{u}(0) = \vc{e} + \sum_{m=1}^{\infty} \vc{B}(m)
(\vc{I} - \vc{G}^m)
(\vc{I} - \vc{A} - \ol{\vc{m}}_{A}\vc{g} )^{-1}\vc{e}
+ \sum_{m=1}^{\infty} {m\vc{B}(m) \over -\sigma}\vc{e},\quad
\label{eqn-u(0)}
\\
\vc{u}(k) =  (\vc{I} - \vc{G}^k)(\vc{I} - \vc{A} - \ol{\vc{m}}_{A}\vc{g} )^{-1}\vc{e}
+ {k \over -\sigma}\vc{e},
\qquad k\in\bbN,
\label{eqn-u(k)}
\end{subnumcases}
and thus
\begin{align}
\lim_{k\to\infty}{\vc{u}(k) \over k}
= {1 \over -\sigma}\vc{e}.
\label{lim-u(k)}
\end{align}
\end{LEM}

\begin{proof}
We first prove (\ref{eqn-u(k)}), which immediately leads to the proof of (\ref{lim-u(k)}).  The matrix generating function $\wh{\vc{G}}(z)$ is the minimal nonnegative solution for the matrix equation $\wh{\mathscr{G}}(z) 
= z\sum_{m=-1}^{\infty}\vc{A}(m) \{ \wh{\mathscr{G}}(z) \}^{m+1}$ (see \cite[Theorems~2.2.1 and 2.2.2]{Neut89}) and thus
\begin{equation}
\wh{\vc{G}}(z)
= z \sum_{m=-1}^{\infty}\vc{A}(m) [\wh{\vc{G}}(z) ]^{m+1}.
\label{eqn:matrx-G}
\end{equation} 
Equation~(\ref{eqn:matrx-G}) is rewritten as
\begin{eqnarray}
\wh{\vc{G}}(z)
&=& 
\left[
\vc{I} - z\sum_{m=0}^{\infty}\vc{A}(m) \{ \wh{\vc{G}}(z) \}^m
\right]^{-1} z\vc{A}(-1).
\label{eqn-wh{G}(z)}  
\end{eqnarray}
Similarly, we have (see \cite[Eq.~(2.4.3)]{Neut89})
\begin{eqnarray}
\wh{\vc{G}}_1(z)
&=& 
\left[
\vc{I} - z\sum_{m=0}^{\infty}\vc{A}(m) \{ \wh{\vc{G}}(z) \}^m
\right]^{-1} z\vc{B}(-1).
\label{eqn-wh{G}_0(z)}
\end{eqnarray}
Combining (\ref{eqn-wh{G}(z)}), (\ref{eqn-wh{G}_0(z)}), and (\ref{cond-finite-MG1-03}) yields
\begin{equation}
\wh{\vc{G}}_1(z)\vc{e}
= \wh{\vc{G}}(z)\vc{e}.
\label{eqn-wh{G}_0(z)e}
\end{equation}
Substituting (\ref{eqn-wh{G}_0(z)e}) into (\ref{eqn-u(k)-02}) and using $\wh{\vc{G}}(1)\vc{e}=\vc{G}\vc{e}=\vc{e}$, we obtain
\begin{align}
\vc{u}(k) 
&= \left. {\rd \over \rd z} [ \wh{\vc{G}}(z) ]^k \right|_{z=1} \vc{e}
= \sum_{n=0}^{k-1} \vc{G}^n \left. {\rd \over \rd z} \wh{\vc{G}}(z) \right|_{z=1} \vc{e} ,
\qquad k \in \bbN.
\label{eqn-u(k)-04}
\end{align}
Note here (see \cite[Eqs.~(3.1.3), (3.1.12), and (3.1.14)]{Neut89}) that
\begin{eqnarray}
\left. {\rd \over \rd z} \wh{\vc{G}}(z) \right|_{z=1} \vc{e}
&=& 
(\vc{I} - \vc{G} + \vc{e}\vc{g})
(\vc{I} - \vc{A} - \ol{\vc{m}}_{A}\vc{g} )^{-1}\vc{e}
\nonumber
\\
&=& 
(\vc{I} - \vc{G} )
(\vc{I} - \vc{A} - \ol{\vc{m}}_{A}\vc{g} )^{-1}\vc{e}
+ {1 \over -\sigma}\vc{e},
\label{eqn-u(1)}
\end{eqnarray}
where the second equality is due to $\vc{g}(\vc{I} - \vc{A} - \ol{\vc{m}}_{A}\vc{g})^{-1}=\vc{\varpi}/(-\sigma)$ (see \cite[Eq.~(3.1.15)]{Neut89}). Inserting (\ref{eqn-u(1)}) into (\ref{eqn-u(k)-04}) results in
\begin{eqnarray*}
\vc{u}(k) 
&=& \sum_{n=0}^{k-1} \vc{G}^n(\vc{I} - \vc{G})
(\vc{I} - \vc{A} - \ol{\vc{m}}_{A}\vc{g} )^{-1}\vc{e}
+ {k \over -\sigma}\vc{e}
\nonumber
\\
&=& (\vc{I} - \vc{G}^k)(\vc{I} - \vc{A} - \ol{\vc{m}}_{A}\vc{g} )^{-1}\vc{e}
+ {k \over -\sigma}\vc{e},
\qquad k \in \bbN,
\end{eqnarray*}
which shows that (\ref{eqn-u(k)}) holds.

Next, we prove (\ref{eqn-u(0)}). Let $\wh{\vc{K}}(z)$, $z \in [0,1]$, denote an $M_0 \times M_0$ matrix such that
\begin{alignat}{2}
(\wh{\vc{K}}(z))_{i,j}
&= \EE_{(0,i)} \!\left [z^{\tau_0}\one(J_{\tau_0} = j) \right],
\qquad i,j\in\bbM_0.
\label{defn-wh{K}(z)}
\end{alignat}
From (\ref{defn-u(k,i)}), we then have
\begin{eqnarray}
\vc{u}(0) 
&=& \left. {\rd \over \rd z} \wh{\vc{K}}(z) \right|_{z=1} \vc{e}.
\label{eqn-u(0)-02}
\end{eqnarray}
We also have (see \cite[Eq (2.4.8)]{Neut89})
\begin{eqnarray}
\wh{\vc{K}}(z)
&=& z\vc{B}(0) 
+ z\sum_{m=1}^{\infty}\vc{B}(m) [ \wh{\vc{G}}(z) ]^{m-1} \wh{\vc{G}}_1(z).
\label{eqn-wh{K}(z)}  
\end{eqnarray}
Combining this and (\ref{eqn-wh{G}_0(z)e}) yields
\begin{align}
\wh{\vc{K}}(z)\vc{e}
&= z\sum_{m=0}^{\infty}\vc{B}(m) [ \wh{\vc{G}}(z) ]^m \vc{e}.
\label{eqn-wh{K}(z)e}  
\end{align}
Substituting (\ref{eqn-wh{K}(z)e}) into (\ref{eqn-u(0)-02}) and using (\ref{eqn-u(k)-04}) and $\sum_{m=0}^{\infty}\vc{B}(m)\vc{e} = \vc{e}$, we obtain
\begin{align}
\vc{u}(0) 
&= \sum_{m=0}^{\infty}\vc{B}(m)\vc{e}
+ \sum_{m=1}^{\infty}\vc{B}(m) 
\left. {\rd \over \rd z} [ \wh{\vc{G}}(z) ]^m \right|_{z=1} \vc{e}
\nonumber
\\
&= \vc{e} + \sum_{m=1}^{\infty}\vc{B}(m) \vc{u}(m).
\label{eqn-u(0)-00}
\end{align}
Finally, inserting (\ref{eqn-u(k)}) into (\ref{eqn-u(0)-00}) leads to 
(\ref{eqn-u(0)}). The proof is completed.
\end{proof}

Lemma~\ref{lem-u(k)}, together with Theorem~\ref{thm-mean}, yields Theorem~\ref{thm-u(k)} below, which ensures that the second moment condition (Assumption~\ref{assum-2nd-moment}) is equivalent to the finiteness of the mean first passage time to level zero in steady state.
\begin{THM}\label{thm-u(k)}
Suppose that Assumption~\ref{assumpt-ergodic-MG1} is satisfied. 
Assumption~\ref{assum-2nd-moment} 
holds if and only if
\begin{equation}
\sum_{k=0}^{\infty} \vc{\pi}(k) \vc{u}(k) < \infty.
\label{sum-k*u(k)e-finite}
\end{equation}
\end{THM}

\begin{proof}
The constant $-\sigma$ in (\ref{lim-u(k)}) is positive and finite. Indeed, it follows from (\ref{defn-beta_A}) and Assumption~\ref{assumpt-ergodic-MG1} that
\begin{align*}
0 
< -\sigma 
&= - \vc{\varpi} \sum_{k=-1}^{\infty} k\vc{A}(k)\vc{e}
\le \vc{\varpi} \vc{A}(-1)\vc{e} 
\le \vc{\varpi} \vc{A}\vc{e} = \vc{\varpi}\vc{e} = 1.
\end{align*}
Therefore, (\ref{sum-k*u(k)e-finite}) is equivalent to (\ref{sum-k*pi(k)e<b}). Furthermore, (\ref{sum-k*pi(k)e<b}) is equivalent to Assumption~\ref{assum-2nd-moment} [due to Theorem~\ref{thm-mean}~(i)]. The proof has been completed. %\qed
\end{proof}

\section{Fundamental deviation matrix of the infinite-level chain}\label{sec:deviation-matrix}

This section consists of two subsections. Section~\ref{subsec:definition-H} defines the fundamental deviation matrix of the infinite-level chain as a block-decomposition-friendly solution for the Poisson equation solved by the deviation matrix (see, e.g., \cite{Cool02}). The fundamental deviation matrix has the same function as the deviation matrix in the perturbation analysis of Markov chains. Furthermore, the fundamental deviation matrix always exists, whereas the deviation matrix does not necessarily, provided that the infinite-level chain is irreducible and positive recurrent. In this sense, we can say that the fundamental deviation matrix is more {\it fundamental} than the deviation matrix. Section~\ref{subsec:decomposition-H} presents an explicit block-decomposition expression of the fundamental deviation matrix. The block-decomposition expression contributes to the analysis of the next section.

\subsection{A block-decomposition-friendly solution for the Poisson equation}\label{subsec:definition-H}

In this subsection, we first define the fundamental deviation matrix after some preparations. We then show that the fundamental deviation matrix is a solution for the Poisson equation of the deviation matrix. That solution, the fundamental deviation matrix, has a finite base set and is suitable for block decomposition by choosing the base set appropriately according to the block structure. We also establish upper bounds for the fundamental deviation matrix and any solution for the Poisson equation. Finally, we show a relationship between the fundamental deviation matrix and the deviation matrix.

We make some preparations to describe the fundamental deviation matrix. 
Let $\bbA$ denote an arbitrary finite subset of $\bbS$, and $\bbB = \bbS \setminus \bbA$, and partition $\vc{P}$ as
\begin{equation}
\vc{P} = \bordermatrix{
  &
\bbA &
\bbB
\cr
\bbA &
\vc{P}_{\bbA} & 
\vc{P}_{\bbA,\bbB}
\cr
\bbB &
\vc{P}_{\bbB,\bbA} & 
\vc{P}_{\bbB}
}.
\end{equation}
Let $\wt{\vc{P}}_{\bbA}$ denote
\begin{align}
\wt{\vc{P}}_{\bbA} = \vc{P}_{\bbA} + \vc{P}_{\bbA,\bbB}(\vc{I} - \vc{P}_{\bbB})^{-1}\vc{P}_{\bbB,\bbA},
\end{align}
which is considered the transition probability matrix of a Markov chain obtained by observing the infinite-level chain $\{(X_n,J_n)\}$ when it is in $\bbA$.
Let $\wt{\vc{\pi}}_{\bbA} :=(\wt{\pi}_{\bbA}(k,i))_{(k,i) \in \bbA}$ denote the stationary distribution vector of $\wt{\vc{P}}_{\bbA}$ and thus
\begin{align*}
\wt{\pi}_{\bbA}(k,i)
= {\pi(k,i) \over \sum_{(\ell,j) \in \bbA}\pi(\ell,j)},
\qquad (k,i) \in \bbA.
\end{align*}
Finally, let $T(\bbA) = \inf\{n \in \bbN: (X_n,J_n) \in \bbA\}$ and $u_{\bbA}(k,i) = \EE_{(k,i)}[T(\bbA)]$ for $(k,i) \in \bbS$.

We are now ready to define the fundamental deviation matrix of the infinite-level chain.
Fix $\vc{\alpha}:=(\alpha_1,\alpha_2) \in \bbS$ arbitrarily, and let $\vc{H}_{\bbA}:=(H_{\bbA}(k,i;\ell,j))_{(k,i;\ell,j) \in \bbS^2}$ denote
\begin{align}
H_{\bbA}(k,i;\ell,j)
&=
\EE_{(k,i)}\!\!
\left[
\sum_{n=0}^{\tau_{\vc{\alpha}}-1} \ol{1}_{(\ell,j)}(X_n, J_n)
\right]
\nonumber
\\
&\qquad {} - \sum_{(k_0,i_0) \in \bbA}
\wt{\pi}_{\bbA}(k_0,i_0) \EE_{(k_0,i_0)}\!\!
\left[
\sum_{n=0}^{\tau_{\vc{\alpha}}-1} \ol{1}_{(\ell,j)}(X_n, J_n)
\right],
\label{defn-H}
\end{align}
where $\ol{1}_{(\ell,j)}(k,i) = 1_{(\ell,j)}(k,i) - \pi(\ell,j)$ for $(k,i;\ell,j) \in \bbS$. Since the infinite-level chain $\{(X_n,J_n)\}$ is irreducible and positive recurrent, we have $\EE_{(k,i)}[ \tau_{\vc{\alpha}}] < \infty$ for any $(k,i) \in \bbS$ and $\vc{\alpha} \in \bbS$. Furthermore, since $\bbA$ is a finite subset of $\bbS$, the second term on the right-hand side of (\ref{defn-H}) is finite and thus $\vc{H}_{\bbA}$ is well-defined. We refer to $\vc{H}_{\bbA}$ as the {\it fundamental deviation matrix with finite base set $\bbA$} or just as the {\it fundamental deviation matrix}.

The fundamental deviation matrix $\vc{H}_{\bbA}$ is closely related to the {\it deviation matrix} $\vc{D}$ defined as follows (though $\vc{D}$ requires some additional conditions; see Proposition~\ref{prop-D-H}):
\begin{align}
\vc{D}
&= \sum_{k=0}^{\infty} ( \vc{P}^k - \vc{e}\vc{\pi} ).
\label{defn:D}
\end{align}
The deviation matrix $\vc{D}$ (if it exists) is a well-known solution for a {\it constrained} Poisson equation (see \cite[Lemma~2.7]{Dend13}):
\begin{align}
(\vc{I} - \vc{P}) \vc{X} &= \vc{I} - \vc{e}\vc{\pi},
\label{Poisson-EQ-H}
\\
\vc{\pi}\vc{X} &=\vc{0}.
\label{Poisson-EQ-H-constraint}
\end{align}
As shown in Theorem~\ref{thm:H}, the fundamental deviation matrix $\vc{H}_{\bbA}$ is a solution for the Poisson equation (\ref{Poisson-EQ-H}) with another constraint:
\begin{align}
( \wt{\vc{\pi}}_{\bbA},~\vc{0} )\vc{X} &=\vc{0}.
\label{Poisson-EQ-HA-constraint}
\end{align}
\begin{THM}\label{thm:H}
Suppose that Assumption~\ref{assumpt-ergodic-MG1} holds, and fix $\vc{\alpha}=(\alpha_1,\alpha_2) \in \bbS$ arbitrarily. For any finite $\bbA \subset \bbS$, $\vc{H}_{\bbA}$ is the unique solution for the Poisson equation (\ref{Poisson-EQ-H}) with the constraint (\ref{Poisson-EQ-HA-constraint}). Moreover, $\vc{H}_{\bbA}$ is independent of $\vc{\alpha}$, and $\vc{H}_{\bbA}(\bbA) := (H_{\bbA}(k,i;\ell,j))_{(k,i;\ell,j) \in \bbA \times \bbS}$ and $\vc{H}_{\bbA}(\bbB) := (H_{\bbA}(k,i;\ell,j))_{(k,i;\ell,j) \in \bbB \times \bbS}$ are given by
\begin{subequations}\label{eqn:H_A}
\begin{align}
\vc{H}_{\bbA}(\bbA) 
&= (\vc{I} - \wt{\vc{P}}_{\bbA} + \vc{e}\wt{\vc{\pi}}_{\bbA})^{-1}
\left[
(\vc{I},~ \vc{P}_{\bbA,\bbB}(\vc{I} - \vc{P}_{\bbB})^{-1})
- \vc{u}_{\bbA}(\bbA)\vc{\pi}
\right],
\label{eqn:H_A(A)}
\\
\vc{H}_{\bbA}(\bbB) 
&= 
(\vc{O},~ (\vc{I} - \vc{P}_{\bbB})^{-1})
- \vc{u}_{\bbB}(\bbA)\vc{\pi}
+ (\vc{I} - \vc{P}_{\bbB})^{-1}\vc{P}_{\bbB,\bbA}\vc{H}_{\bbA}(\bbA),
\label{eqn:H_A(B)}
\end{align}
\end{subequations}
where $\vc{u}_{\bbA}(\bbA) = (u_{\bbA}(k,i))_{(k,i) \in \bbA}$ and $\vc{u}_{\bbA}(\bbB) = (u_{\bbA}(k,i))_{(k,i) \in \bbB}$.
\end{THM}

\begin{REM}
Theorem~\ref{thm:H} is immediately extended to the general class of Markov chains with countable states because the proof of this theorem does not depend on the M/G/1-type structure. Related results are found in \cite[Theorem~3.2]{Jinpeng_Liu22}. 
\end{REM}

\begin{REM}
The matrix $\vc{H}_{\bbA}$ can be expressed explicitly in a block-decomposition-form with probabilistically interpretable matrices, given that the base set $\bbA$ is appropriately chosen (see Corollary~\ref{coro:dev-MG1}). 
\end{REM}

{\it Proof of Theorem~\ref{thm:H}.}
We first prove that $\vc{H}_{\bbA}$ is a solution for the Poisson equation (\ref{Poisson-EQ-H}) with (\ref{Poisson-EQ-HA-constraint}). To facilitate our discussion, we define $\wt{\vc{H}}:=(\wt{H}(k,i;\ell,j))_{(k,i;\ell,j) \in \bbS^2}$ as
\begin{align}
\wt{H}(k,i;\ell,j)
&=\EE_{(k,i)}\!\!
\left[
\sum_{n=0}^{\tau_{\vc{\alpha}}-1} \ol{1}_{(\ell,j)}(X_n, J_n)
\right]
\nonumber
\\
&=
\EE_{(k,i)}\!\!
\left[
\sum_{n=0}^{\tau_{\vc{\alpha}}-1} 1_{(\ell,j)}(X_n, J_n)
\right]
- \pi(\ell,j)\EE_{(k,i)}[\tau_{\vc{\alpha}}].
\label{defn-wt{H}}
\end{align}
We also define
\begin{align*}
\wt{\vc{H}}(\bbA) = (\wt{H}(k,i;\ell,j))_{(k,i;\ell,j) \in \bbA \times \bbS},
\qquad
\wt{\vc{H}}(\bbB) = (\wt{H}(k,i;\ell,j))_{(k,i;\ell,j) \in \bbB \times \bbS}.
\end{align*}
It then follows from (\ref{defn-H}) and (\ref{defn-wt{H}}) that
\begin{align}
\vc{H}_{\bbA}
= \wt{\vc{H}} 
- \vc{e} \cdot \wt{\vc{\pi}}_{\bbA}\wt{\vc{H}}(\bbA).
\label{add:eqn:H_A-wt{H}}
\end{align}
Solving the matrix equation (\ref{add:eqn:H_A-wt{H}}) for $\wt{\vc{H}}$ and 
decomposing its rows with the sets $\bbA$ and $\bbB$, we have
\begin{subequations}\label{eqn:H_A-wt{H}}
\begin{align}
\wt{\vc{H}}(\bbA) 
&=  \vc{H}_{\bbA}(\bbA) + \vc{e}\wt{\vc{\pi}}_{\bbA}\wt{\vc{H}}(\bbA),
\label{eqn:H_A-wt{H}-a}
\\
\wt{\vc{H}}(\bbB)
&=  \vc{H}_{\bbA}(\bbB) + \vc{e}\wt{\vc{\pi}}_{\bbA}\wt{\vc{H}}(\bbA),
\label{eqn:H_A-wt{H}-b}
\end{align}
\end{subequations}
where the sizes of the vectors $\vc{e}$ in (\ref{eqn:H_A-wt{H}-a}) and (\ref{eqn:H_A-wt{H}-b}) are equal to the cardinalities of $\bbA$ and $\bbB$, respectively. Note (see \cite[Lemma~2.1]{Dend13}) that $\wt{\vc{H}}$ is a solution for the Poisson equation (\ref{Poisson-EQ-H}) such that $\wt{H}(\vc{\alpha};\ell,j)=0$ for all $(\ell,j) \in \bbS$. Therefore $\vc{H}_{\bbA}$, satisfying (\ref{add:eqn:H_A-wt{H}}), is a solution for the Poisson equation (\ref{Poisson-EQ-H}). In addition, it follows from (\ref{add:eqn:H_A-wt{H}}) and $\wt{\vc{\pi}}_{\bbA}\vc{e}=1$ that
\begin{align}
( \wt{\vc{\pi}}_{\bbA},~\vc{0} )\vc{H}_{\bbA} =\vc{0},
\label{cond:wt{pi}_A}
\end{align}
and thus 
$\vc{H}_{\bbA}$ satisfies the constraint (\ref{Poisson-EQ-HA-constraint}).

Next, we prove that (\ref{eqn:H_A(A)}) and (\ref{eqn:H_A(B)}) hold. Since the matrix $\wt{\vc{H}}$, defined in (\ref{defn-wt{H}}), is a solution for the Poisson equation (\ref{Poisson-EQ-H}), it follows from \cite[Theorem~2.5]{Dend13}) that
\begin{subequations}\label{add:eqn:H_A}
\begin{align}
(\vc{I} - \wt{\vc{P}} _{\bbA})
\wt{\vc{H}}(\bbA) &= 
(\vc{I},~ \vc{P}_{\bbA,\bbB}(\vc{I} - \vc{P}_{\bbB})^{-1})
- \vc{u}_{\bbA}(\bbA)\vc{\pi},
\label{add:eqn:H_A(A)}
\\
\wt{\vc{H}}(\bbB) &= 
(\vc{O},~ (\vc{I} - \vc{P}_{\bbB})^{-1})
- \vc{u}_{\bbB}(\bbA)\vc{\pi}
+ (\vc{I} - \vc{P}_{\bbB})^{-1}\vc{P}_{\bbB,\bbA}\wt{\vc{H}}(\bbA).
\label{add:eqn:H_A(B)}
\end{align}
\end{subequations}
Substituting (\ref{eqn:H_A-wt{H}-a}) into (\ref{add:eqn:H_A(A)}) yields
\begin{align}
(\vc{I} - \wt{\vc{P}} _{\bbA})
\vc{H}_{\bbA}(\bbA) &= 
(\vc{I},~ \vc{P}_{\bbA,\bbB}(\vc{I} - \vc{P}_{\bbB})^{-1})
- \vc{u}_{\bbA}(\bbA)\vc{\pi}.
\label{cond:wh{H}(A)-01}
\end{align}
Furthermore, (\ref{cond:wt{pi}_A}) yields
\begin{align}
\wt{\vc{\pi}}_{\bbA}\vc{H}_{\bbA}(\bbA) &= \vc{0}.
\label{cond:wh{H}(A)-02}
\end{align}
Note here that $\vc{I} - \wt{\vc{P}}_{\bbA} + \vc{e}\wt{\vc{\pi}}_{\bbA}$ is non-singular and
\begin{align}
(\vc{I} - \wt{\vc{P}}_{\bbA} + \vc{e}\wt{\vc{\pi}}_{\bbA})^{-1}
(\vc{I} - \wt{\vc{P}}_{\bbA})
= \vc{I} - \vc{e}\wt{\vc{\pi}}_{\bbA}.
\label{eqn:220508-01}
\end{align}
Therefore, pre-multiplying by $(\vc{I} - \wt{\vc{P}}_{\bbA} + \vc{e}\wt{\vc{\pi}}_{\bbA})^{-1}$ the both sides of (\ref{cond:wh{H}(A)-01}) and using (\ref{cond:wh{H}(A)-02}) and (\ref{eqn:220508-01}), we obtain
\begin{align*}
\vc{H}_{\bbA}(\bbA) &= 
(\vc{I} - \wt{\vc{P}}_{\bbA} + \vc{e}\wt{\vc{\pi}}_{\bbA})^{-1}
\left[
(\vc{I},~ \vc{P}_{\bbA,\bbB}(\vc{I} - \vc{P}_{\bbB})^{-1})
- \vc{u}_{\bbA}(\bbA)\vc{\pi}
\right],
\end{align*}
which shows that (\ref{eqn:H_A(A)}) holds. In addition, substituting (\ref{eqn:H_A-wt{H}}) into (\ref{add:eqn:H_A(B)}) and using $(\vc{I} - \vc{P}_{\bbB})^{-1}\vc{P}_{\bbB,\bbA}\vc{e}=\vc{e}$, we obtain
\begin{align*}
&
\vc{H}_{\bbA}(\bbB) + \vc{e}\wt{\vc{\pi}}_{\bbA}\wt{\vc{H}}(\bbA) 
\nonumber
\\
&\quad = 
(\vc{O},~ (\vc{I} - \vc{P}_{\bbB})^{-1})
- \vc{u}_{\bbB}(\bbA)\vc{\pi}
+ (\vc{I} - \vc{P}_{\bbB})^{-1}\vc{P}_{\bbB,\bbA}
\left[
\vc{H}_{\bbA}(\bbA) + \vc{e}\wt{\vc{\pi}}_{\bbA}\wt{\vc{H}}(\bbA)
\right]
\nonumber
\\
&\quad = 
(\vc{O},~ (\vc{I} - \vc{P}_{\bbB})^{-1})
- \vc{u}_{\bbB}(\bbA)\vc{\pi}
+ (\vc{I} - \vc{P}_{\bbB})^{-1}\vc{P}_{\bbB,\bbA}
\vc{H}_{\bbA}(\bbA)
+ \vc{e}\wt{\vc{\pi}}_{\bbA}\wt{\vc{H}}(\bbA),
\end{align*}
and thus
\begin{align*}
\vc{H}_{\bbA}(\bbB)
= 
(\vc{O},~ (\vc{I} - \vc{P}_{\bbB})^{-1})
- \vc{u}_{\bbB}(\bbA)\vc{\pi}
+ (\vc{I} - \vc{P}_{\bbB})^{-1}\vc{P}_{\bbB,\bbA}
\vc{H}_{\bbA}(\bbA),
\end{align*}
which shows that (\ref{eqn:H_A(B)}) holds.

Equation (\ref{eqn:H_A}) implies that $\vc{H}_{\bbA}$ is independent of $\vc{\alpha}$. To complete the proof, it thus suffices to show that $\vc{H}_{\bbA}$ is the unique solution for the Poisson equation (\ref{Poisson-EQ-H}) with (\ref{Poisson-EQ-HA-constraint}). Based on (\ref{cond:wh{H}(A)-01}) and (\ref{cond:wh{H}(A)-02}), we can consider $\vc{H}_{\bbA}(\bbA)$ to be a solution for
\begin{subequations}\label{eqn:wh{H}(A)-02}
\begin{align}
(\vc{I} - \wt{\vc{P}} _{\bbA})
\vc{X}(\bbA) &= 
(\vc{I},~ \vc{P}_{\bbA,\bbB}(\vc{I} - \vc{P}_{\bbB})^{-1})
- \vc{u}_{\bbA}(\bbA)\vc{\pi},
\\
\wt{\vc{\pi}}_{\bbA}\vc{X}(\bbA) &= \vc{0}.
\end{align}
\end{subequations}
We now define $\vc{X}_1(\bbA)$ and $\vc{X}_2(\bbA)$ as arbitrary solutions for (\ref{eqn:wh{H}(A)-02}), and then define $\vc{h}_1(\bbA)$ and $\vc{h}_2(\bbA)$ as arbitrary columns of $\vc{X}_1(\bbA)$ and $\vc{X}_2(\bbA)$, respectively. By definition, $\vc{h}(\bbA):=\vc{h}_1(\bbA)-\vc{h}_2(\bbA)$ is harmonic for $\wt{\vc{P}} _{\bbA}$, i.e., $\vc{h}(\bbA) = \wt{\vc{P}} _{\bbA}\vc{h}(\bbA)$, which yields
\begin{align*}
\vc{h}(\bbA) = {1 \over n} \sum_{k=1}^n \wt{\vc{P}} _{\bbA}^k \vc{h}(\bbA),
\qquad n \in \bbN.
\end{align*}
Combining this with \cite[Theorem~14.3.6]{Meyn09} leads to $\vc{h}(\bbA) = \wt{\vc{\pi}}_{\bbA}\vc{h}(\bbA) \cdot \vc{e}$ (see also \cite[Proposition~17.4.1]{Meyn09}) and thus
\begin{align*}
\vc{X}_1(\bbA) - \vc{X}_2(\bbA) = \vc{e} \vc{\eta}(\bbA)
\quad \mbox{for some row vector $\vc{\eta}(\bbA)$}.
\end{align*}
Furthermore, by definition, $\wt{\vc{\pi}}_{\bbA}(\vc{X}_1(\bbA) - \vc{X}_2(\bbA)) = \vc{0}$ and thus $\vc{\eta}(\bbA) = \vc{0}$, which implies $\vc{X}_1(\bbA) - \vc{X}_2(\bbA) = \vc{O}$. Therefore, $\vc{H}_{\bbA}(\bbA)$ is the unique solution for the set of equations (\ref{eqn:wh{H}(A)-02}), and $\vc{H}_{\bbA}(\bbB)$ is uniquely determined by (\ref{eqn:H_A(B)}). Consequently, $\vc{H}_{\bbA}$ is the unique solution for the Poisson equation (\ref{Poisson-EQ-H}) with (\ref{Poisson-EQ-HA-constraint}). The proof has been completed.

\smallskip

The following lemma presents respective upper bounds for $|\vc{H}_{\bbA}|\vc{e}$ and any solution of the Poisson equation (\ref{Poisson-EQ-H}) through a drift condition different from the one given in Lemma~\ref{lem-drift-cond}.
\begin{LEM}\label{lem:Pv'}
Suppose that Assumption~\ref{assumpt-ergodic-MG1} is satisfied. Let $\vc{v}':=(v'(k,i))_{(k,i) \in \bbS}$ denote a column vector such that
\begin{equation}
\vc{v}'(k):= (v'(k,i))_{i\in\bbM_{k\vmin1}}
= 
\left\{
\begin{array}{ll}
\dm{1 \over -\sigma}\vc{e},				& k=0,
\\
\rule{0mm}{6mm}
\dm{1 \over -\sigma}(k\vc{e} + \vc{a}), 	& k \in\bbN,
\end{array}
\right.
\label{defn-v'(k)}
\end{equation}
where $\vc{a} \ge \vc{0}$ is given in (\ref{defn-alpha}). The following then hold.
\begin{enumerate}
\item There exist some $b' \in (0,\infty)$ and $K' \in \bbN$ such that
\begin{alignat}{2}
\vc{P}^{(N)}\vc{v}' \le \vc{P}\vc{v}' &\le \vc{v}' - \vc{e} + b' \vc{1}_{\bbL_{\le K'}}& \quad &\mbox{for all $N \in \bbN$}.
\label{drift-cond-02a}
\end{alignat}
\item Any solution $\vc{X}$ for the Poisson equation (\ref{Poisson-EQ-H}) satisfies $|\vc{X}| \le C_0\vc{v}'\vc{e}^{\top}$ for some $C_0 > 0$. Furthermore, $\vc{\pi}|\vc{X}| < \infty$ under Assumption~\ref{assum-2nd-moment}. 
\item For each finite $\bbA \subset \bbS$, there exists some $C_{\bbA} > 0$ such that $|\vc{H}_{\bbA}|\vc{e} \le  C_{\bbA}\vc{v}'$.
\end{enumerate}
\end{LEM}

\begin{proof}
See Appendix~\ref{proof:lem:Pv'}.
\end{proof}

As mentioned above, the following proposition shows that $\vc{D}$ does not necessarily exist even though $\vc{H}_{\bbA}$ does. 
\begin{PROP}\label{prop-D-H}
Suppose that Assumption~\ref{assumpt-ergodic-MG1} is satisfied and $\vc{P}$ is aperiodic. 
\begin{enumerate}
\item Assumption~\ref{assum-2nd-moment} holds if and only if $\vc{D}$ exists.
\item If $\vc{D}$ exists, then it is the unique solution (\ref{Poisson-EQ-H}) with (\ref{Poisson-EQ-H-constraint}) such that $\vc{\pi} |\vc{X}| < \infty$, and furthermore, for any finite $\bbA \subset \bbS$,
\begin{align}
\vc{D} = (\vc{I} - \vc{e}\vc{\pi})\vc{H}_{\bbA}.
\label{eqn:D}
\end{align}
\end{enumerate}
\end{PROP}

\begin{proof}
See Appendix~\ref{appen:proof-prop-D-H}.
%\qed
\end{proof}

\begin{REM}
A similar and general version of Proposition~\ref{prop-D-H} is presented in \cite[Corollary~3.1]{Jinpeng_Liu22}.
\end{REM}

Although the deviation matrix has been used for the perturbation analysis of block-structured Markov chains (see, e.g., \cite{Dend13, LiuYuan14, LiuYuan18-STM}), we take the fundamental deviation matrix $\vc{H}_{\bbA}$ instead of the deviation matrix $\vc{D}$. This is because $\vc{D}$ requires the aperiodicity of $\vc{P}$ but our analysis does not necessarily require it.

\subsection{Explicit block decomposition of the fundamental deviation matrix}\label{subsec:decomposition-H}

The fundamental deviation matrix $\vc{H}_{\bbA}$ with $\bbA=\bbL_0$ has an explicit block- decomposition form with probabilistically interpretable matrices. Thus $\vc{H}_{\bbL_0}$ facilitates the analysis developed in the next section. For simplicity, we omit the subscript ``$\bbL_0$" of $\vc{H}_{\bbL_0}$ and define $\vc{H}(k;\ell)$ as
\begin{alignat*}{2}
\vc{H}(k;\ell)
&=(H(k,i;\ell,j))_{(k,i;\ell,j) \in \bbM_{k\vmin1} \times \bbM_{\ell\vmin1}}
&\qquad k,\ell \in \bbZ_+.
\end{alignat*}

The following is a corollary of Theorem~\ref{thm:H}, which expresses $\vc{H}$ in a block-decomposition-form with probabilistically interpretable matrices.
\begin{COR}\label{coro:dev-MG1}
Suppose that Assumption~\ref{assumpt-ergodic-MG1} holds. We then have
%
%\begin{subequations}
\begin{alignat}{2}
\vc{H}(0;\ell)  
&= 
(\vc{I} - \vc{K} + \vc{e}\vc{\kappa})^{-1}
(\vc{I} - \vc{u}(0) \vc{\pi}(0))\vc{F}_+(0;\ell), &  \ell &\in \bbZ_+,
\label{eqn-H(0;l)}
\\
\vc{H}(k;\ell)
&= (1 - \delta_{0,\ell}) \vc{F}_{+}(k;\ell) - \vc{u}(k) \vc{\pi}(\ell)
\nonumber
\\
& \qquad  {} + \vc{G}^{k-1} (\vc{I}-\vc{\Phi}(0))^{-1} \vc{B}(-1)\vc{H}(0;\ell),
& \quad k \in \bbN, ~\ell &\in \bbZ_+,
\label{eqn-H(k;l)}
\end{alignat}
%\end{subequations}
%
where $\vc{\kappa}$ denotes the unique stationary distribution vector of $\vc{K}:=\wh{\vc{K}}(1)$, and where 
$\vc{F}_+(k;\ell)$, $k,\ell\in\bbZ_+$, denotes an $\bbM_{k \vmin 1} \times \bbM_{\ell \vmin 1}$ matrix such that
\begin{equation*}
(\vc{F}_+(k;\ell))_{i,j}
= 
\EE_{(k,i)}\!\!
\left[ 
\sum_{n=0}^{\tau_0-1} 1_{(\ell,j)}(X_n, J_n)
\right],\qquad k,\ell\in\bbZ_+,
%\label{eqn:F_+^{(N)}(k;l)}
\end{equation*}
and thus, for all $k,\ell \in \bbZ_+$ such that $k \vmin \ell = 0$,
\begin{subnumcases}{\vc{F}_+(k;\ell)=\label{eqn-F_+}}
\vc{I}, & $k = 0,~\ell = 0$,\label{eqn-wt{H}-a}
\\
\dm\sum_{m=1}^{\infty} \vc{B}(m)\vc{F}_{+}(m;\ell), 
& $k=0,~\ell \in \bbN$,\label{eqn-wt{H}-b}
\\
\vc{O}, & $k \in \bbN,~\ell=0$.
\label{eqn-wt{H}-c}
\end{subnumcases}
\end{COR}

\begin{proof}
We first prove (\ref{eqn-H(0;l)}). By definition, $\vc{K} = \wt{\vc{P}}_{\bbL_0}$ and thus $\vc{\kappa}$ is the stationary probability vector of $\wt{\vc{P}}_{\bbL_0}$. Therefore, letting $\bbA = \bbL_0$ in Theorem~\ref{thm:H}, we have 
\begin{align}
\vc{H}(0;\ell)  
&= (\vc{I} - \vc{K} + \vc{e}\vc{\kappa})^{-1} 
[\vc{F}_+(0;\ell) - \vc{u}(0)\vc{\pi}(\ell)],
\label{eqn-wt{H}(0;l)}
\\
\vc{H}(k;\ell) &= \vc{F}_+(k;\ell) - \vc{u}(k)\vc{\pi}(\ell)
+ \vc{F}_+(k;1)\vc{B}(-1)\vc{H}(0;\ell),
\qquad k \in \bbN.
\label{eqn-wt{H}(k;l)}
\end{align}
It follows from (\ref{eqn-F_+}) and \cite[Theorem~10.0.1]{Meyn09}) that
\begin{align}
\vc{\pi}(\ell)
= \vc{\pi}(0)\vc{F}_+(0;\ell),
\qquad \ell \in \bbZ_+.
\end{align}
Substituting this into (\ref{eqn-wt{H}(0;l)}) yields (\ref{eqn-H(0;l)}).

Next, we prove (\ref{eqn-H(k;l)}). Note (\cite[Theorem 9]{Zhao03}) that 
\begin{equation}
\vc{F}_{+}(k;\ell)
= \vc{G}^{k-\ell}\vc{F}_{+}(\ell;\ell), 
\qquad k \in \bbZ_{\ge \ell},~\ell \in\bbN,
\label{eq-H(m;k)}
\end{equation}
and thus 
\begin{equation}
\vc{F}_{+}(k;1)
= \vc{G}^{k-1}\vc{F}_{+}(1;1)
=  \vc{G}^{k-1}(\vc{I} - \vc{\Phi}(0))^{-1},
\qquad k \in \bbN,
\label{eqn-F_{+}(k;1)}
\end{equation}
where $\vc{\Phi}(0)$ is given in (\ref{defn-Phi(0)}). Note also that $\vc{F}_{+}(k;0) = \vc{O}$ for $k\in\bbN$, which is shown in (\ref{eqn-wt{H}-c}). To emphasize this exceptional case, we write
\begin{align}
\vc{F}_{+}(k;\ell) = (1 - \delta_{0,\ell})\vc{F}_{+}(k;\ell),
\qquad k \in \bbN,~\ell\in\bbZ_+.
\label{add-eqn:F(k;l)}
\end{align}
Substituting (\ref{eqn-F_{+}(k;1)}) and (\ref{add-eqn:F(k;l)}) into (\ref{eqn-wt{H}(k;l)}), we obtain (\ref{eqn-H(k;l)}). The proof has been completed. %\qed
\end{proof}

\begin{rem}
Let $\wt{\vc{H}}(k;\ell)
=(\wt{H}(k,i;\ell,j))_{(k,i;\ell,j) \in \bbM_{k\vmin1} \times \bbM_{\ell\vmin1}}$ for $k,\ell \in \bbZ_+$. Letting $\bbA = \bbL_0$ in (\ref{add:eqn:H_A}) and following the proof of Corollary~\ref{coro:dev-MG1}, we obtain the equations similar to (\ref{eqn-H(0;l)}) and (\ref{eqn-H(k;l)}):
\begin{alignat*}{2}
(\vc{I} - \vc{K})\wt{\vc{H}}(0;\ell)  
&= 
(\vc{I} - \vc{u}(0) \vc{\pi}(0))\vc{F}_+(0;\ell), &  \ell &\in \bbZ_+,
\\
\wt{\vc{H}}(k;\ell)
&= (1 - \delta_{0,\ell}) \vc{F}_{+}(k;\ell) - \vc{u}(k) \vc{\pi}(\ell)
\nonumber
\\
& \qquad  {} + \vc{G}^{k-1} (\vc{I}-\vc{\Phi}(0))^{-1} \vc{B}(-1)\wt{\vc{H}}(0;\ell),
& \quad k \in \bbN, ~\ell &\in \bbZ_+.
\end{alignat*}
\end{rem}

%%%%%%%%%%%%%%%%%%%%%%%%% Convergence %%%%%%%%%%%%%%%%%%%%%%%%%%%
 
\section{Subgeometric convergence in the infinite-level limit}\label{sec-convergence}

The main purpose of this section is to prove the subgeometric convergence of the level-wise difference $\vc{\pi}^{(N)}(k) - \vc{\pi}(k)$. This section consists of three subsections. Section~\ref{subsec:difference-formula} presents a difference formula for $\vc{\pi}^{(N)}$ and $\vc{\pi}$, and then Section~\ref{subsec:basic-convergence} shows the uniform convergence of $\|\vc{\pi}^{(N)} - \vc{\pi}\|$ under Assumption~\ref{assum-2nd-moment}. Based on those results, Section~\ref{sec-subgeometric} derives a subgeometric convergence formula for $\vc{\pi}^{(N)}(k) - \vc{\pi}(k)$ under an additional condition (Assumption~\ref{assumpt-tail}). The subgeometric convergence formula is presented in Theorem~\ref{thm-asym}, which is the main theorem of this paper.

\subsection{A difference formula of the finite- and infinite-level stationary distributions}\label{subsec:difference-formula}

This subsection presents a difference formula for $\vc{\pi}^{(N)}$ and $\vc{\pi}$ with the fundamental deviation matrix $\vc{H}$. The following two matrices are needed to describe the difference formula.
\begin{align}
\ool{A}(k) 
&= \sum_{\ell=k+1}^{\infty} \ol{\vc{A}}(\ell), 
\quad
\ool{B}(k) = \sum_{\ell=k+1}^{\infty} \ol{\vc{B}}(\ell), 
\quad 
k \in \bbZ_{\ge -1}.
\label{defn:ool{A}(n)-ool{B}(n)}
\end{align} 
where $\ol{\vc{A}}(k)$ and $\ol{\vc{B}}(k)$ are given in (\ref{defn-ol{A}(k)}) and (\ref{defn-ol{B}(k)}), respectively.
\begin{LEM}[Difference formula]\label{prop-diff}
Suppose that Assumption~\ref{assumpt-ergodic-MG1} holds, then
\begin{align}
\vc{\pi}^{(N)} - \vc{\pi} 
= \vc{\pi}^{(N)} ( \vc{P}^{(N)} - \vc{P} ) \vc{H},
\qquad N \in \bbN.
\label{diff-formula}
\end{align}
For $k \in \bbZ_{[0,N]}$, we also have
\begin{align}
&
\vc{\pi}^{(N)}(k)-\vc{\pi}(k)
\nonumber
\\
&\quad=
\vc{\pi}^{(N)}(0) 
\left[
\sum_{n=N+1}^{\infty} \vc{B}(n) \vc{S}^{(N)}(n;k)
  + {1 \over -\sigma} \ool{B}(N-1) \vc{e}\vc{\pi}(k)
\right] 
\notag\\
&\quad\quad
+ \sum_{\ell=1}^N \vc{\pi}^{(N)}(\ell) 
\left[
\sum_{n=N+1}^{\infty} \vc{A}(n-\ell) \vc{S}^{(N)}(n;k)
+ {1 \over -\sigma} \ool{A}(N-\ell-1) \vc{e}\vc{\pi}(k)
\right],
\label{eq-xN(k)-x(k)}
\end{align}
where
\begin{align}
\vc{S}^{(N)}(n;k) 
&=
(1 - \delta_{0,k})(\vc{G}^{N-k}-\vc{G}^{n-k}) \vc{F}_{+}(k;k)
\nonumber
\\
& {}
+ (\vc{G}^{N-1}-\vc{G}^{n-1})(\vc{I}-\vc{\Phi}(0))^{-1} \vc{B}(-1)\vc{H}(0;k) 
\nonumber
\\
& {}
+ (\vc{G}^{N}-\vc{G}^{n}) 
\left( \vc{I}-\vc{A} - \ol{\vc{m}}_{A}\vc{g} \right)^{-1} \vc{e}
 \vc{\pi}(k),\quad n \in \bbZ_{\ge N+1},~k \in \bbZ_{[0,N]}.
\label{def-S(N;n,k)}
\end{align}
\end{LEM}

\begin{proof}
First, we prove (\ref{diff-formula}). From (\ref{Poisson-EQ-H}), $\vc{\pi}^{(N)} =\vc{\pi}^{(N)}\vc{P}^{(N)}$, and $\vc{\pi}\vc{e}=1$, we obtain
\begin{eqnarray*}
\vc{\pi}^{(N)}  ( \vc{P}^{(N)} - \vc{P} ) \vc{X}
&=& \vc{\pi}^{(N)} ( \vc{I} - \vc{P} ) \vc{X}
= \vc{\pi}^{(N)} ( \vc{I} - \vc{e}\vc{\pi} )
= \vc{\pi}^{(N)} - \vc{\pi},
\end{eqnarray*}
which shows that (\ref{diff-formula}) holds. 

Next, we prove (\ref{eq-xN(k)-x(k)}). To this end, we begin with expressing the difference formula (\ref{diff-formula}) level-wise. From (\ref{defn-P}) and (\ref{defn:P^{(N)}}), we  have
\begin{eqnarray*}
\lefteqn{
\vc{P}^{(N)} - \vc{P} 
}
\quad &&
\nonumber\\
&=&
\!
\bordermatrix{
& \bigcup_{\ell=0}^{N-1} \bbL_{\ell}
& \bbL_N  			
& \bbL_{N+1} 		
& \bbL_{N+2}		
& \cdots
\cr
\bbL_0			&
\vc{O}      	& 
\ol{\vc{B}}(N) &  
 -\vc{B}(N+1)  	& 
-\vc{B}(N+2) 	& 
\cdots 
\cr
\bbL_1			&
\vc{O}    		& 
\ol{\vc{A}}(N-1) & 
-\vc{A}(N)  	&  
-\vc{A}(N+1)	& 
\cdots 			
\cr
\vdots 			& 
\vdots 			& 
\vdots 			& 
\vdots 			&
\vdots 			
\cr
\bbL_{N-1}\! 	&
\vc{O}    		& 
\ol{\vc{A}}(1)  & 
-\vc{A}(2)   	&    
-\vc{A}(3)		& 
\cdots 
\cr
\bbL_N			&
\vc{O}     		& 
\ol{\vc{A}}(0)  & 
-\vc{A}(1)      & 
-\vc{A}(2)      & 
\cdots 
\cr
\bbL_{N+1}\! 	&
\vc{O}     		& 
-\vc{A}(-1)  	& 
-\vc{A}(0)   	& 
-\vc{A}(1)   	& 
\cdots 
\cr
\bbL_{N+2}\! 	&
\vc{O}  		& 
\vc{O}  		& 
-\vc{A}(-1)   	& 
-\vc{A}(0)   	& 
\cdots 
\cr
~~~\vdots      	& 
\vdots         	& 
\vdots          & 
\vdots          & 
\vdots          & 
\ddots 
}.
\label{hm-eqn-01}
\qquad
\end{eqnarray*}
Using this equation, we decompose (\ref{diff-formula}) into level-wise expressions: 
For $k \in \bbZ_{[0,N]}$,
\begin{eqnarray}
\vc{\pi}^{(N)}(k) - \vc{\pi}(k) 
&=&
\left[
\vc{\pi}^{(N)}(0) \ol{\vc{B}}(N)
+ \sum_{\ell=1}^N \vc{\pi}^{(N)}(\ell)\ol{\vc{A}}(N-\ell)
\right]
\vc{H}(N;k) 
\nonumber
\\
&&{} 
-\sum_{n=N+1}^{\infty} 
\left[
\vc{\pi}^{(N)}(0)\vc{B}(n) + \sum_{\ell=1}^N \vc{\pi}^{(N)}(\ell) \vc{A}(n-\ell) 
\right]
 \vc{H}(n;k)
\nonumber
\\
&=& 
\vc{\pi}^{(N)}(0) 
\sum_{n=N+1}^{\infty} \vc{B}(n) 
\left[ \vc{H}(N;k) - \vc{H}(n;k) \right]
\nonumber
\\
&&{} +
\sum_{\ell=1}^N \vc{\pi}^{(N)}(\ell)
\sum_{n=N+1}^{\infty} \vc{A}(n-\ell) 
\left[ \vc{H}(N;k)-\vc{H}(n;k) \right].\qquad
\label{eq-kankei}
\end{eqnarray}

To proceed further, we rewrite the term $\vc{H}(N;k) - \vc{H}(n;k)$ in (\ref{eq-kankei}) by using $\vc{G}$ and related matrices (including vectors) introduced in Section~\ref{sec-second-order-moment}. Combining (\ref{eqn-H(k;l)}) and (\ref{eq-H(m;k)}), we obtain,
for $n \in \bbZ_{\ge N+1}$ and $k \in \bbZ_{[0,N]}$, 
\begin{eqnarray}
\vc{H}(N;k)-\vc{H}(n;k) 
&=& 
(1 - \delta_{0,k})(\vc{G}^{N-k}-\vc{G}^{n-k}) \vc{F}_{+}(k;k)
\nonumber
\\
&& {}
+ (\vc{G}^{N-1}-\vc{G}^{n-1}) (\vc{I}-\vc{\Phi}(0))^{-1} \vc{B}(-1)\vc{H}(0;k) 
\nonumber
\\
&& {}
+ \left[ \vc{u}(n)-\vc{u}(N) \right] \vc{\pi}(k).
\label{eq-K-K-tmp}
\end{eqnarray}
Using (\ref{eqn-u(k)}), we rewrite $\vc{u}(n)-\vc{u}(N)$ in (\ref{eq-K-K-tmp}) as
\begin{align}
\vc{u}(n)-\vc{u}(N)
&= (\vc{G}^{N}-\vc{G}^{n}) 
(\vc{I} - \vc{A} - \ol{\vc{m}}_{A}\vc{g} )^{-1} \vc{e}
+ {1 \over -\sigma}(n-N)\vc{e},~~~ n \in \bbZ_{\ge N+1}.
\label{eq-mu(N)-mu(n+ell)}
\end{align}
Applying (\ref{eq-mu(N)-mu(n+ell)}) to (\ref{eq-K-K-tmp}) and using (\ref{def-S(N;n,k)}), we obtain, for $n \in \bbZ_{\ge N+1}$ and $k \in \bbZ_{[0,N]}$, 
\begin{eqnarray}
\vc{H}(N;k)-\vc{H}(n;k) 
&=& 
(1 - \delta_{0,k})(\vc{G}^{N-k}-\vc{G}^{n-k}) \vc{F}_{+}(k;k)
\nonumber
\\
&& {}
+ (\vc{G}^{N-1}-\vc{G}^{n-1}) (\vc{I}-\vc{\Phi}(0))^{-1} \vc{B}(-1)\vc{H}(0;k)
\nonumber
\\ 
&& {}
+ (\vc{G}^{N}-\vc{G}^{n}) 
(\vc{I} - \vc{A} - \ol{\vc{m}}_{A}\vc{g} )^{-1} \vc{e}\vc{\pi}(k)
\nonumber
\\
&& {}
+ {1 \over -\sigma} (n-N) \vc{e}\vc{\pi}(k)
\nonumber
\\ 
&=& \vc{S}^{(N)}(n;k) + {1 \over -\sigma} (n-N) \vc{e}\vc{\pi}(k).
\label{eq-K-K}
\end{eqnarray}
Furthermore, substituting (\ref{eq-K-K}) into (\ref{eq-kankei}) results in the following: For $k \in \bbZ_{[0,N]}$,
\begin{eqnarray}
\lefteqn{
\vc{\pi}^{(N)}(k) - \vc{\pi}(k) 
}
~~&&
\nonumber
\\
&=& 
\vc{\pi}^{(N)}(0) 
\sum_{n=N+1}^{\infty} \vc{B}(n) 
\left\{ \vc{S}^{(N)}(n;k) + {1 \over -\sigma} (n-N) \vc{e}\vc{\pi}(k)\right\}
\nonumber
\\
&&{} +
\sum_{\ell=1}^N \vc{\pi}^{(N)}(\ell)
\sum_{n=N+1}^{\infty} \vc{A}(n-\ell) 
\left\{ \vc{S}^{(N)}(n;k) + {1 \over -\sigma} (n-N) \vc{e}\vc{\pi}(k)\right\}.\qquad
\label{eq-kankei-tmp}
\end{eqnarray}

To prove (\ref{eq-xN(k)-x(k)}), we arrange the terms of (\ref{eq-kankei-tmp}). For $\ell\in\bbZ_{[1,N]}$, we have
\begin{align}
&
\sum_{n=N+1}^\infty (n-N) \vc{A}(n-\ell)
\nonumber
\\
&\quad=
\sum_{m=1}^{\infty} m \vc{A}(m + N - \ell)
= 
\sum_{m=1}^{\infty} \sum_{k=1}^m \vc{A}(m + N - \ell)
\nonumber
\\
&\quad= \sum_{k=1}^{\infty} \sum_{m=k}^{\infty} \vc{A}(m + N - \ell)
= \sum_{k=1}^{\infty} \ol{\vc{A}}(k + N - \ell - 1)
\nonumber
\\
&\quad= \ool{\vc{A}}(N - \ell - 1),
\label{eq-(n-N)A}
\end{align}
and, similarly, 
\begin{align} 
\sum_{n=N+1}^{\infty} (n-N)\vc{B}(n) &= \ool{B}(N-1).
\label{eq-(n-N)B}
\end{align}
Substituting (\ref{eq-(n-N)A}) and (\ref{eq-(n-N)B}) into (\ref{eq-kankei-tmp}) yields (\ref{eq-xN(k)-x(k)}). The proof is completed.
\end{proof}

\begin{REM}\label{rem:deviation-matrix}
In the same way as for (\ref{diff-formula}), we can derive a similar difference formula by any solution $\vc{X}$ for the Poisson equation (\ref{Poisson-EQ-H}):
\begin{equation}
\vc{\pi}^{(N)} - \vc{\pi}
= \vc{\pi}^{(N)} ( \vc{P}^{(N)} - \vc{P} )\vc{X},
\qquad N \in \bbN.
\label{diff-formula-D^{(N)}}
\end{equation}
Note that (\ref{diff-formula-D^{(N)}}) holds for $\vc{X}=\wt{\vc{H}}$ defined in (\ref{defn-wt{H}}) and for $\vc{X}=\vc{D}$ defined in (\ref{defn:D})  (if exists); the latter case is presented in \cite[Equation~(3)]{Heid08}.
\end{REM}

\subsection{Uniform convergence under the second-order moment condition}\label{subsec:basic-convergence}

This subsection presents two theorems. The first theorem shows 
the uniform convergence of $\{\vc{\pi}^{(N)}\}$ to $\vc{\pi}$ under the second-order moment condition (Assumption~\ref{assum-2nd-moment}). 
The second one shows a useful inequality for the tail probabilities of $\vc{\pi}=(\vc{\pi}(0),\vc{\pi}(1),\dots)$ and $\vc{\pi}^{(N)}=(\vc{\pi}^{(N)}(0),\vc{\pi}^{(N)}(1),\dots,\vc{\pi}^{(N)}(N))$.

The following is the first theorem.
\begin{THM}\label{thm-asymp-pi^(N)(k)}
If Assumptions~\ref{assumpt-ergodic-MG1} and \ref{assum-2nd-moment} hold, then $\| \vc{\pi}^{(N)} - \vc{\pi} \| = o(N^{-1})$.
\end{THM}

\begin{proof}
It suffices to show that
\begin{align}
\sum_{k=0}^N |\vc{\pi}^{(N)}(k)-\vc{\pi}(k)|\vc{e}
= o(N^{-1}).
\label{asymp:sum-pi^{(N)}(k)-pi(k)}
\end{align}
This is because 
\begin{align*}
\| \vc{\pi}^{(N)} - \vc{\pi} \| 
&= \sum_{k=0}^N |\vc{\pi}^{(N)}(k)-\vc{\pi}(k)|\vc{e}
+ \sum_{k=N+1}^{\infty} \vc{\pi}(k)\vc{e},
\end{align*}
and $\vc{\pi}(k) =\ul{\bcalO}(k^{-2})$ [due to Theorem~\ref{thm-mean}~(i)]. 

In what follows, we prove (\ref{asymp:sum-pi^{(N)}(k)-pi(k)}). 
Equation (\ref{defn-v'(k)}) shows that $\vc{v}'(n)> \vc{v}'(N)$ for $n \in\bbZ_{\ge N+1}$, and thus Lemma~\ref{lem:Pv'}~(iii) yields
\begin{align}
\sum_{k=0}^N \big[ |\vc{H}(N;k)|\vc{e} + |\vc{H}(n;k)|\vc{e} \big]
\le 2C\vc{v}'(n),\qquad n \in\bbZ_{\ge N+1},
\label{eqn:210309-01}
\end{align}
where $C > 0$ is some constant. It follows from (\ref{eqn:210309-01}) and (\ref{eq-kankei}) that
\begin{align}
&
\sum_{k=0}^N |\vc{\pi}^{(N)}(k)-\vc{\pi}(k)|\vc{e}
\nonumber
\\
&\quad\le
2C
\left[
\vc{\pi}^{(N)}(0) 
\sum_{n=N+1}^{\infty} \vc{B}(n) \vc{v}'(n)
+
\sum_{\ell=1}^N \vc{\pi}^{(N)}(\ell)
\sum_{n=N+1}^{\infty} \vc{A}(n-\ell) \vc{v}'(n)
\right].
\label{eqm:200524-01}
\end{align}
It also follows from (\ref{defn-v'(k)}), $\vc{A}(k) = \ul{\bcalO}(k^{-3})$, and $\vc{B}(k) = \ul{\bcalO}(k^{-3})$ (see Remark~\ref{rem:assum-2nd-moment}) that
\begin{align}
\sum_{n=N+1}^{\infty} \vc{B}(n) \vc{v}'(n)
&= {1 \over -\sigma}
\sum_{n=N+1}^{\infty} \vc{B}(n) (n\vc{e} + \vc{a})
= \ul{\bcalO}(N^{-1}),
\label{asymp:sum-nB(n)}
\\
\sum_{n=N+1}^{\infty} \vc{A}(n-\ell) \vc{v}'(n)
&= {1 \over -\sigma}
\sum_{n=N+1}^{\infty} \vc{A}(n-\ell) (n\vc{e} + \vc{a})
\nonumber
\\
&= {1 \over -\sigma}
\sum_{n=N+1}^{\infty} \ul{\bcalO}((n-\ell)^{-3}) 
[ (n-\ell)\vc{e} + \ell\vc{e} + \vc{a} ]
\nonumber
\\
&= 
\ul{\bcalO}((N-\ell+1)^{-1})
+ \ell  \ul{\bcalO}((N-\ell+1)^{-2}),~~ \ell \in \bbZ_{[1,N]}.
\label{asymp:sum-(n-l)A(n-l)}
\end{align}
Estimating the right-hand side of (\ref{eqm:200524-01}) by (\ref{asymp:sum-nB(n)}), (\ref{asymp:sum-(n-l)A(n-l)}), and $\vc{\pi}^{(N)}(\ell) = \ul{\bcalO}(\ell^{-2})$ (due to Theorem~\ref{thm-mean}~(ii)), we obtain
\begin{align*}
&
\sum_{k=0}^N |\vc{\pi}^{(N)}(k)-\vc{\pi}(k)|\vc{e}
\nonumber
\\
&\quad \le
o(N^{-1})
+
\sum_{\ell=1}^N 
o(\ell^{-2})
\left[
o((N-\ell+1)^{-1})
+ \ell \cdot o((N-\ell+1)^{-2})
\right]
\nonumber
\\
&\quad=  
o(N^{-1})
+ \sum_{\ell=1}^N 
o(\ell^{-2}) o((N-\ell+1)^{-1})
+ \sum_{\ell=1}^N 
o(\ell^{-1}) o((N-\ell+1)^{-2})
\nonumber
\\
&\quad= o(N^{-1}),
\end{align*}
which shows that (\ref{asymp:sum-pi^{(N)}(k)-pi(k)}) holds.
The proof is completed. 
\end{proof}

The following theorem is the second one in this subsection.
\begin{THM}\label{THM:bound-ol{pi}^{(N)}(k)}
Suppose that Assumptions~\ref{assumpt-ergodic-MG1} and \ref{assum-2nd-moment} are satisfied. For $k \in \bbZ_+$, let
 $\ol{\vc{\pi}}(k):=(\ol{\pi}(k,i))_{i \in \bbM_1}$ and $\ol{\vc{\pi}}^{(N)}(k):=(\ol{\pi}^{(N)}(k,i))_{i \in \bbM_1}$ denote
\begin{align}
\ol{\vc{\pi}}(k) 
&=\sum_{\ell=k+1}^{\infty}\vc{\pi}(\ell),  \qquad
\ol{\vc{\pi}}^{(N)}(k) 
=\sum_{\ell=k+1}^{\infty}\vc{\pi}^{(N)}(\ell),
\label{defn:ol{pi}(k)-pi^{(N)}(k)}
\end{align}
respectively, where $\vc{\pi}^{(N)}(k) = \vc{0}$ for $k \in \bbZ_{\ge N+1}$. We then have
\begin{equation}
(\ol{\vc{\pi}}^{(N)}(0), \ol{\vc{\pi}}^{(N)}(1),\ol{\vc{\pi}}^{(N)}(2), \dots)
\le (1+ o(N^{-1})) (\ol{\vc{\pi}}(0), \ol{\vc{\pi}}(1), \ol{\vc{\pi}}(2), \dots).
\label{asymp-pi^{(N)}(k)}
\end{equation}
\end{THM}

\begin{proof}
See Appendix~\ref{proof:THM:bound-ol{pi}^{(N)}(k)}.
\end{proof}

\subsection{Subgeometric convergence}\label{sec-subgeometric}

In this subsection, we first provide preliminaries on heavy-tailed distributions and subgeometric functions, and then present the main theorem and its corollaries on the subgeometric convergence of $\{\vc{\pi}^{(N)}\}$.

\subsubsection{Preliminaries: Heavy-tailed distributions and subgeometric functions}\mbox{}

We introduce the class of heavy-tailed distributions on the domain $\bbZ_+$ together with and its two subclasses. To this end, let $F:\bbZ_+\to [0,1]$ denote a probability distribution (function), and let $\ol{F} = 1 - F$, which is the complementary distribution (function) of $F$. Furthermore, let $F^{*2}(k) = \sum_{\ell=0}^k F(\ell) F(k-\ell)$ for $k \in \bbZ_+$.
\begin{DEF}[Heavy-tailed, long-tailed, and subexponential distributions]
\hfill
\begin{enumerate}
\item  A distribution $F$ is said to be heavy-tailed if and only if
\begin{equation*}
\limsup_{k \to \infty} e^{\theta k} \ol{F}(k) = \infty
\quad \mbox{for any $\theta > 0$}.
%\label{defn:class-H}
\end{equation*}
The class of heavy-tailed distributions is denoted by $\calH$.
\item A distribution $F$ is said to be {\it long-tailed} if and only if
\begin{equation*}
\lim_{k \to \infty} 
{\ol{F}(k+\ell) \over \ol{F}(k)} = 1
\quad \mbox{for any fixed $\ell \in \bbN$}.
%\label{defn:class-L}
\end{equation*}
The class of long-tailed distributions is denoted by $\calL$.
\item A distribution $F$ is said to be {\it subexponential} if and only if
\begin{equation*}
\lim_{k \to \infty} 
{1 - F^{*2}(k) \over \ol{F}(k)} = 2.
\end{equation*}
The class of  subexponential distributions is denoted by $\calS$.
\end{enumerate}
\end{DEF}

\begin{REM}\label{rem:classes-H_L_S}
The inclusion relation of the above three classes holds: $\calS \subsetneq \calL \subsetneq \calH$ (see \cite[Lemmas 2.17 and 3.2; Section~3.7]{Foss13}).
\end{REM}

Next, we introduce the class of subgeometric functions. 
\begin{DEF}[Subgeometric functions]\label{defn-subgeo}
A function $r:\bbZ_+\to\bbR_+:=[0,\infty)$ is said to be a {\it subgeometric function} if and only if $\log r(k) = o(k)$ as $k \to \infty$.
The class of subgeometric functions is denoted by $\varTheta$.
\end{DEF}

\begin{PROP}\label{prop-subgeo}
If $F \in \calL$, then $\ol{F} \in \varTheta$.
\end{PROP}

\begin{proof}
It suffices to show that
\begin{align}
-\log \ol{F}(k) = o(k)\quad \mbox{as $k\to\infty$}.
\label{eqn:210311-01}
\end{align}
Let $R(k) = -\log \ol{F}(k)$ for $k \in \bbZ_+$. It then follows (see \cite[Lemma~2.22]{Foss13}) that $\lim_{k\to\infty} [R(k+1) - R(k)] = 0$. Therefore, for any $\varep > 0$, there exists some $k_0 \in \bbN$ such that $|R(k+1) - R(k)| < \varep$ for all $k \in \bbZ_{\ge k_0}$, which yields
\begin{align*}
{R(k_0) - \varep n \over k_0+n} 
<
{R(k_0+n) \over k_0+n} < {R(k_0) + \varep n \over k_0+n}\quad \mbox{for all $n \in \bbN$}.
\end{align*}
Letting $n\to\infty$ and then $\varep \downarrow 0$ in the above inequality, we obtain
\begin{align*}
\lim_{n\to\infty}{R(k_0+n) \over k_0+n} = 0,
\end{align*}
which implies that (\ref{eqn:210311-01}) holds. The proof is completed.
\end{proof}

\begin{REM}
In relation to class $\varTheta$, there is a class $\varLambda$ of {\it subgeometric rate functions} introduced in \cite{Numm83}. A function $\psi:\bbZ_+\to\bbR_+$ belongs to class $\varLambda$ if and only if 
\[
0 <
\liminf_{k\to\infty}{\psi(k) \over \psi_0(k)} \le
\limsup_{k\to\infty}{\psi(k) \over \psi_0(k)} < \infty,
\]
for some $\psi_0:\bbZ_+\to [2,\infty)$ such that $\psi_0$ is nondecreasing such that $\log \psi_0(k)/k \searrow 0$ as $k \to \infty$. By definition, $\varLambda \subset \varTheta$.
\end{REM}

\subsubsection{The main theorem and its corollaries}\mbox{}

We now make an additional assumption on $\{\vc{A}(k)\}$ and $\{\vc{B}(k)\}$ to study the subgeometric convergence of $\{\vc{\pi}^{(N)}\}$ to $\vc{\pi}$. 
\begin{ASSU} \label{assumpt-tail}
There exists a distribution $F \in \calS$ such that 
\begin{subequations}\label{eq-assum-tail}
\begin{eqnarray}
\lim_{k\to\infty} { \ool{\vc{A}}(k)\vc{e} \over \ol{F}(k)} 
&=& \vc{c}_{A}, 
\label{asymp-ool{A}(k)}
\\
\lim_{k\to\infty} {\ool{\vc{B}}(k)\vc{e} \over \ol{F}(k)} 
&=& \vc{c}_{B}, 
\label{asymp-ool{B}(k)}
\end{eqnarray}
\end{subequations}
where either $\vc{c}_{A}\neq\vc{0}$ or $\vc{c}_{B}\neq\vc{0}$.
\end{ASSU}

Assumption~\ref{assumpt-tail} can be interpreted as a condition on the distribution of level increments in steady state. We define $\Delta_+ = \max(X_1-X_0,0)$ and call it the {\it nonnegative level increment}. We then define $D(k)$, $k \in \bbZ_+$, as
\begin{align}
D(k)
&= \sum_{(\ell,i) \in \bbS} \pi(\ell,i) 
\PP(\Delta_+ \le k \mid (X_0,J_0)=(\ell,i))
\nonumber
\\
&= \sum_{(\ell,i) \in \bbS} \pi(\ell,i) 
\PP(X_1-X_0 \le k \mid (X_0,J_0)=(\ell,i))
\nonumber
\\
&= \vc{\pi}(0) \sum_{n=0}^k \vc{B}(n)\vc{e}
+ \ol{\vc{\pi}}(0) \sum_{n=-1}^k \vc{A}(n)\vc{e}.
\label{defn:D(k)}
\end{align}
We call $D$ the {\it stationary nonnegative level-increment (SNL) distribution}. Assumption~\ref{assumpt-ergodic-MG1} ensures that the SNL distribution $D$ has a finite positive mean $\mu_{\pi}:=\vc{\pi}(0) \ol{\vc{m}}_{B} + \ol{\vc{\pi}}(0) \ol{\vc{m}}_{A}^+ > 0$, where 
\begin{align*}
\ol{\vc{m}}_{A}^+ 
= \sum_{k=1}^{\infty}k\vc{A}(k)\vc{e}
= \ol{\vc{m}}_{A} + \vc{A}(-1)\vc{e}
\ge \vc{0},\neq\vc{0}.
\end{align*}
We also define $\ol{D}_I(k) = 1 - D_I(k)$ for $k \in \bbZ_+$, where $D_I$ denotes the integrated tail distribution (the equilibrium distribution) of the SNL distribution $D$, that is,
\begin{align}
D_I(k)
= \mu_{\pi}^{-1} \sum_{\ell=0}^k (1 - D(\ell)),
\qquad k \in \bbZ_+.
\label{defn:D_I(k)}
\end{align}
It follows from (\ref{defn:D_I(k)}), (\ref{defn:D(k)}), and Assumption \ref{assumpt-tail} that 
\begin{align}
\lim_{k\to\infty}{\ol{D}_I(k) \over \ol{F}(k)}
= {
\vc{\pi}(0) \vc{c}_{B} + \ol{\vc{\pi}}(0) \vc{c}_{A}
\over
\vc{\pi}(0) \ol{\vc{m}}_{B} + \ol{\vc{\pi}}(0) \ol{\vc{m}}_{A}^+
} \in (0,\infty).
\label{asymp:F_pi^I}
\end{align}
Since $F \in \calS$,  we have $D_I \in \calS \subset \calL$ (see \cite[Corollary~3.13]{Foss13}) and thus $\ol{D}_I \in \varTheta$ due to Proposition~\ref{prop-subgeo}.

Assumption~\ref{assumpt-tail} yields a subexponential asymptotic formula for the stationary distribution $\vc{\pi}$ of the infinite-level chain, as shown in the next proposition. This proposition contributes to the proof of Theorem~\ref{thm-asym} below (see (\ref{eqn:210216-06b}) in Appendix~\ref{proof:lem:main-01}).
\begin{PROP}[{}{\cite[Theorem~3.1]{Masu16-ANOR}}] \label{prop-Masu16}
If Assumptions~\ref{assumpt-ergodic-MG1} and \ref{assumpt-tail} hold, then
\begin{equation}
\lim_{k\to\infty} {\ol{\vc{\pi}}(k) \over \ol{F}(k)}
= { 
\vc{\pi}(0) \vc{c}_{B} + \ol{\vc{\pi}}(0) \vc{c}_{A}
\over
-\sigma 
} \vc{\varpi}.
\label{asymp:ol{pi}(n)}
\end{equation}
\end{PROP}

We have arrived at the main theorem of this paper.
\begin{THM} \label{thm-asym}
If Assumptions~\ref{assumpt-ergodic-MG1}, \ref{assum-2nd-moment}, and \ref{assumpt-tail} hold, then
\begin{alignat}{2}
\lim_{N\to\infty}
{
\vc{\pi}^{(N)}(k) - \vc{\pi}(k)  
\over
\ol{F}(N)
}
&=
{
\vc{\pi}(0)\vc{c}_{B} + \ol{\vc{\pi}}(0)\vc{c}_{A}
\over -\sigma
} \vc{\pi}(k), &\qquad k &\in \bbZ_+,
\label{main_thm_forluma-01}
\end{alignat}
and thus, as $N\to\infty$, $\vc{\pi}^{(N)}(k) - \vc{\pi}(k)$ converge to zero according to a subgeometric function $\ol{F}$ (see Proposition~\ref{prop-subgeo}).
\end{THM}

\begin{proof}
See Appendix~\ref{proof:thm-asym}.
\end{proof}

There are two other versions of the subgeometric convergence formula (\ref{main_thm_forluma-01}).
\begin{COR}\label{cor:subgeo_convergence_01}
If all the conditions of Theorem~\ref{thm-asym} are satisfied, then
\begin{alignat}{2}
\lim_{N\to\infty}
{
\vc{\pi}^{(N)}(k) - \vc{\pi}(k)  
\over
\ol{D}_I(N) 
}
&=
{\vc{\pi}(0) \ol{\vc{m}}_{B} + \ol{\vc{\pi}}(0) \ol{\vc{m}}_{A}^+
 \over -\sigma}
\vc{\pi}(k), &\qquad k &\in \bbZ_+,
\label{main_thm_forluma-02}
\\
\lim_{N\to\infty}
{
\vc{\pi}^{(N)}(k) - \vc{\pi}(k)  
\over
\ol{\vc{\pi}}(N)\vc{e}
}
&= \vc{\pi}(k), &\qquad k &\in \bbZ_+.
\label{main_thm_forluma-03}
\end{alignat}
\end{COR}

\begin{proof}
The formulas (\ref{main_thm_forluma-02}) and (\ref{main_thm_forluma-03}) follow from the combination of Theorem~\ref{thm-asym} with (\ref{asymp:F_pi^I}) and (\ref{asymp:ol{pi}(n)}), respectively.
\end{proof}

Additionally, Theorem~\ref{thm-asym} (together with Corollary~\ref{cor:subgeo_convergence_01}) yields three types of subgeometric convergence formulas for the relative difference $[\vc{\pi}^{(N)}(k) - \vc{\pi}(k)] / \vc{\pi}(k)\vc{e}$.
\begin{COR}\label{cor:subgeo_convergence_02}
If all the conditions of Theorem~\ref{thm-asym} are satisfied, then
\begin{alignat*}{2}
\lim_{N\to\infty}
{1 \over 
\ol{F}(N)
}
{
\left\|
\vc{\pi}^{(N)}(k) - \vc{\pi}(k)  
\over
\vc{\pi}(k)\vc{e}
\right\|
}
&=
{
\vc{\pi}(0)\vc{c}_{B} + \ol{\vc{\pi}}(0)\vc{c}_{A}
\over -\sigma
}, &\qquad k &\in \bbZ_+,
%\label{main_cor_forluma-01}
\\
\lim_{N\to\infty}
{1 \over 
\ol{D}_I(N)
}
{
\left\|
\vc{\pi}^{(N)}(k) - \vc{\pi}(k)  
\over
\vc{\pi}(k)\vc{e}
\right\|
}
&=
{\vc{\pi}(0) \ol{\vc{m}}_{B} + \ol{\vc{\pi}}(0) \ol{\vc{m}}_{A}^+
 \over -\sigma}, &\qquad k &\in \bbZ_+,
%\label{main_cor_forluma-02}
\\
\lim_{N\to\infty}
{1 \over \ol{\vc{\pi}}(N)\vc{e}}
{
\left\|
\vc{\pi}^{(N)}(k) - \vc{\pi}(k)  
\over
\vc{\pi}(k)\vc{e}
\right\|
}
&= 1, &\qquad k &\in \bbZ_+.
%\label{main_cor_forluma-03}
\end{alignat*}
\end{COR}

\section{Concluding remarks}\label{sec:remarks}

Theorem~\ref{thm-asym}, the main theorem of this paper, presents the subgeometric convergence formula for the difference between $\vc{\pi}^{(N)}(k)$ and $\vc{\pi}(k)$, which are the respective stationary probabilities of level $k$ in the finite- and infinite-level M/G/1-type Markov chains. 
Theorem~\ref{thm-asym}, together with Corollaries~\ref{cor:subgeo_convergence_01} and \ref{cor:subgeo_convergence_02}, provides three pieces of knowledge on the subgeometric convergence of $\{\vc{\pi}^{(N)}\}$ to $\vc{\pi}$:
\begin{enumerate}
\item The convergence of $\{\vc{\pi}^{(N)}\}$ to $\vc{\pi}$ is subgeometric if the integrated tail distribution $D_I$ of the SNL distribution $D$ is subexponential.
\item The subgeometric convergence speed of $\{\vc{\pi}^{(N)}\}$ is equal to the decay speed of $\ol{D}_I(N)$ and $\ol{\vc{\pi}}(N)$. 
\item The decay speed of $\|\vc{\pi}^{(N)}(k) - \vc{\pi}(k)\| / \vc{\pi}(k)\vc{e}$ is independent of the level value $k$.

\end{enumerate}

There are two challenging problems related to this study. The first problem is to derive geometric convergence formulas for $\vc{\pi}^{(N)}(k) - \vc{\pi}(k)$. One of those formulas would take a form:
\[
\lim_{N\to\infty}
{
\vc{\pi}^{(N)}(k) - \vc{\pi}(k)  
\over
\gamma^N
}
=\vc{\phi}(k), \qquad k \in \bbZ_+,
\]
for some $\gamma \in (0,1)$ and nonnegative vector $\vc{\phi}(k) \neq \vc{0}$. If such a geometric convergence formula is obtained, we can derive a geometric asymptotic formula for the loss probability in M/G/1-type queues. The second problem is to study the convergence speed of the total variation distance $\|\vc{\pi}^{(N)} - \vc{\pi}\|=\sum_{k=0}^{\infty}|\vc{\pi}^{(N)}(k) - \vc{\pi}(k)|\vc{e}$ in the geometric and subgeometric convergence cases. The first problem would be more challenging than the first one because it is not, in general, allowed to interchanging the order of two operators ``$\lim_{N\to\infty}$" and ``$\sum_{k=0}^{\infty}$".

%%%%%%%%%%%%%%%%%%%%%%%%%%%%%%%%%%%%%%%%%%%%%%%%%%%%%%%%%%%%%%%%%%%%%%%%%%%%
%\fi

%%%%%%%%%%%%%%%%%%%%%%%%%%% Appendix %%%%%%%%%%%%%%%%%%%%%%%%%%%%%%
\appendix

\section{Proof of Lemma~\ref{lem:Pv'}}\label{proof:lem:Pv'}

We prove the statement (i). It follows from (\ref{defn-P}), (\ref{defn:P^{(N)}}), and (\ref{defn-v'(k)}) that $\vc{P}^{(N)}\vc{v}' \le \vc{P}\vc{v}'$ for $N \in \bbN$ and that, for $k \in \bbZ_{\ge 2}$,
\begin{eqnarray*}
\sum_{\ell=0}^{\infty}\vc{P}(k;\ell)\vc{v}'(\ell)
&=& {1 \over -\sigma}
\left(
\sum_{\ell=-1}^{\infty} (\ell+k) \vc{A}(\ell)\vc{e}
+ \sum_{\ell=-1}^{\infty} \vc{A}(\ell)\vc{a}
\right)
\nonumber
\\
&=& {1 \over -\sigma}
\left[
k\vc{e} 
+ \left( 
\vc{A}\vc{a} + \vcm_{A}
\right)
\right]
\nonumber
\\
&=& {1 \over -\sigma} (k\vc{e} + \vc{a})
+ {1 \over -\sigma}
\left( 
\vc{A}\vc{a} + \vcm_{A} - \vc{a}
\right)
\nonumber
\\
&=& {1 \over -\sigma} (k\vc{e} + \vc{a}) + {\sigma \over -\sigma}\vc{e}
\nonumber
\\
&=& \vc{v}'(k) - \vc{e},
\end{eqnarray*}
where the second last equality is due to (\ref{eqn:a}).
In a similar way, we can show that $\sum_{\ell=0}^{\infty}\vc{P}(0;\ell)\vc{v}'(\ell) < \infty$ and $\sum_{\ell=0}^{\infty}\vc{P}(1;\ell)\vc{v}'(\ell) < \infty$. 
Therefore, (\ref{drift-cond-02a}) holds for some $b' \in (0,\infty)$ and $K' \in \bbN$.

Next, we prove the statement (ii). Let $\vc{x}_{(\ell,j)}$, $(\ell,j) \in \bbS$, denote the $(\ell,j)$-th column of an arbitrary solution $\vc{X}$ for the Poisson equation (\ref{Poisson-EQ-H}). The $\vc{x}_{(\ell,j)}$ is a solution for a Poisson equation
\begin{align*}
(\vc{I} - \vc{P}) \vc{x}_{(\ell,j)} 
= \vc{1}_{(\ell,j)} - \vc{e}\pi(\ell,j).
\end{align*}
Therefore, it follows from the statement (i) (which has been proved) and \cite[Theorem~2.3]{Glyn96} that there exists some $C_0 > 0$ such that $|\vc{x}_{(\ell,j)}| \le C_0 \vc{v}'$ for any $(\ell,j) \in \bbS$ and thus $|\vc{X}| \le C_0 \vc{v}'\vc{e}^{\top}$.  In addition, it follows from (\ref{defn-v'(k)}) and Theorem~\ref{thm-mean} that if Assumption~\ref{assum-2nd-moment} holds then $\vc{\pi}\vc{v}' < \infty$ and thus $\vc{\pi}|\vc{X}| \le C_0 (\vc{\pi}\vc{v}')\vc{e}^{\top}< \infty$. Consequently, the statement (ii) is true.

In what follows, we prove the statement (iii). To this end, we first confirm that the statement~(iii) holds  if there exists some $C > 0$ such that
\begin{align}
\EE_{(k,i)}\!\!\left[ 
\sum_{n=0}^{\tau_{\vc{\alpha}} - 1}
1_{\bbL_{\le K'}}(X_n, J_n)  \right]
\le C \quad \mbox{for all $(k,i) \in \bbS$}.
\label{eqn:201227-01}
\end{align}
From (\ref{eqn:H_A-wt{H}}), we have
\begin{align}
|\vc{H}_{\bbA}| 
\le |\wt{\vc{H}}| 
+ \vc{e}\left(\wt{\vc{\pi}}_{\bbA} \mid \vc{0} \right)  |\wt{\vc{H}}|. 
\label{bound:|H_A|}
\end{align}
It follows from (\ref{defn-wt{H}}) that
\begin{align}
\sum_{(\ell,j) \in \bbS}
|\wt{H}(k,i;\ell,j)|
&\le 2 \EE_{(k,i)}\!\! \left[ \tau_{\vc{\alpha}}  \right],
\qquad (k,i) \in \bbS.
\label{eqn:201227-02}
\end{align}
It also follows from the statement (i) and \cite[Theorem~2.1]{Glyn96} that
\begin{align}
\EE_{(k,i)}\!\! \left[ \tau_{\vc{\alpha}}  \right]
\le v'(k,i) + b'
\EE_{(k,i)}\!\!\left[ 
\sum_{n=0}^{\tau_{\vc{\alpha}} - 1}
1_{\bbL_{\le K'}}(X_n, J_n)  \right], 
\quad (k,i) \in \bbS.
\label{bound:H^{(N)}(k,i;l,j)}
\end{align}
In addition, (\ref{defn-v'(k)}) yields $\vc{e} \le (-\sigma)\vc{v}'$.
Therefore, if (\ref{eqn:201227-01}) holds, then combining (\ref{eqn:201227-02}) and (\ref{bound:H^{(N)}(k,i;l,j)}) leads to
\begin{align*}
\sum_{(\ell,j) \in \bbS}
|\wt{H}(k,i;\ell,j)|
&\le 2 [v'(k,i) + b'C] 
\le 2 [1 + (-\sigma)b'C] v'(k,i),
\quad (k,i) \in \bbS,
\end{align*}
and thus $|\wt{\vc{H}}|\vc{e} \le C' \vc{v}'$ for some $C' > 0$. Applying this inequality to (\ref{bound:|H_A|}), we have
\begin{align*}
|\vc{H}_{\bbA}| \vc{e} 
&\le C' \vc{v}' 
+ C'\{ \left(\wt{\vc{\pi}}_{\bbA} \mid \vc{0} \right)\vc{v}' \}\vc{e} 
\nonumber
\\
&\le C' \vc{v}' 
+ C' { \vc{\pi}\vc{v}' \over \sum_{(k,i) \in \bbA} \pi(k,i)} \vc{e} 
\nonumber
\\
&\le C'
\left[ 
1 + { (-\sigma) \vc{\pi}\vc{v}' \over \sum_{(k,i) \in \bbA} \pi(k,i)}
\right]
\vc{v}',
\end{align*}
which implies that the statement (iii) holds.

To show (\ref{eqn:201227-01}), we establish a nontrivial inequality for $1_{\bbL_{\le K'}}(k,i)$. Let $\vc{\varGamma}:=(\varGamma(k,i;\ell,j))_{(k,i;\ell,j) \in \bbS^2}$ denote a stochastic matrix such that
\begin{align*}
\varGamma(k,i;\ell,j)
= \sum_{m=0}^{\infty} (1 - \theta)\theta^m P^m(k,i;\ell,j),
\qquad (k,i;\ell,j) \in \bbS^2,
\end{align*}
where $\theta \in (0,1)$ is some constant and $P^m(k,i;\ell,j)$ denotes the $(k,i;\ell,j)$-th element of $\vc{P}^m$. Let $\tau_{\ge k} = \inf\{n\in\bbN:X_n \ge k\}$, and let $N' > K'$ be a positive integer. It then follows that, for $N' > K'$,
\begin{align}
\varGamma(k,i;\vc{\alpha})
&= \EE_{(k,i)}\!\!\left[ 
\sum_{m=0}^{\infty}(1 - \theta)\theta^m
1_{\vc{\alpha}}(X_m, J_m)
\right]
\nonumber
\\
&\ge
\EE_{(k,i)}\!\!\left[ 
\sum_{m=0}^{\tau_{\ge N'} - 1}
(1 - \theta)\theta^m
1_{\vc{\alpha}}(X_m, J_m)
\right], \qquad (k,i) \in \bbL_{\le K'}.
\label{eqn:211014-02}
\end{align}
The irreducibility of $\{(X_n,J_n)\}$ implies that, for each $(k,i) \in \bbL_{\le K'}$, the right-hand side of (\ref{eqn:211014-02}) is nondecreasing with $N'$ and is positive for all sufficiently large $N'> K'$. Therefore, for all sufficiently large $N'>K'$, we have
\begin{align*}
0 < 
\varphi &:= 
\min_{(k,i) \in \bbL_{\le K'}} 
\EE_{(k,i)}\!\!\left[ 
\sum_{m=0}^{\tau_{\ge N'} - 1}
(1 - \theta)\theta^m
1_{\vc{\alpha}}(X_m, J_m)
\right]
\le \min_{(k,i) \in \bbL_{\le K'}} \varGamma(k,i;\vc{\alpha}),
\end{align*}
which leads to the desired inequality for $1_{\bbL_{\le K'}}(k,i)$: 
\begin{align}
1_{\bbL_{\le K'}}(k,i) 
&\le \varphi^{-1}\varGamma(k,i;\vc{\alpha})
= \varphi^{-1}\sum_{m=0}^{\infty}
(1 - \theta)\theta^m P^m(k,i;\vc{\alpha}),\quad (k,i) \in \bbS.
\label{eqn:201018-02}
\end{align} 

Using the inequality (\ref{eqn:201018-02}), we complete the proof. It follows from (\ref{eqn:201018-02}) that
\begin{align*}
\EE_{(k,i)}\!\!\left[ 
\sum_{n=0}^{\tau_{\vc{\alpha}} - 1}
1_{\bbL_{\le K'}}(X_n, J_n)  \right]
&\le 
\EE_{(k,i)}\!\!\left[ 
\sum_{n=0}^{\tau_{\vc{\alpha}} - 1}
\varphi^{-1}
\sum_{m=0}^{\infty} (1 - \theta)\theta^m
P^m(X_n, J_n;\vc{\alpha})
\right]
\nonumber
\\
&= 
\varphi^{-1}
\sum_{m=0}^{\infty} (1 - \theta)\theta^m
\EE_{(k,i)}\!\!\left[ 
\sum_{n=0}^{\tau_{\vc{\alpha}} - 1}
P^m(X_n, J_n;\vc{\alpha})
\right]
\nonumber
\\
&= 
\varphi^{-1}
\sum_{m=0}^{\infty} (1 - \theta)\theta^m
\EE_{(k,i)}\!\!\left[ 
\sum_{n=0}^{\tau_{\vc{\alpha}} - 1}
1_{\vc{\alpha}}(X_{n+m}, J_{n+m}) 
\right]
\nonumber
\\
&\le 
\varphi^{-1}
\sum_{m=0}^{\infty} (1 - \theta)\theta^m m
= {\theta \over \varphi(1 - \theta)},
\qquad (k,i) \in \bbS.
%\label{bound:H^{(N)}(k,i;l,j)-02}
\end{align*}
Consequently, (\ref{eqn:201227-01}) holds for $C = \theta/[\varphi(1 - \theta)] > 0$. The proof is completed.

\section{Proof of Proposition~\ref{prop-D-H}}\label{appen:proof-prop-D-H}

We first prove the statement (ii), assuming that the statement (i) is true. According to this assumption, the existence of $\vc{D}$ implies that Assumption~\ref{assum-2nd-moment} holds. Lemma~\ref{lem:Pv'}~(ii) thus ensures that any solution $\vc{X}$ for the Poisson equation (\ref{Poisson-EQ-H}) satisfies $\vc{\pi} |\vc{X}| < \infty$, and especially, $\vc{\pi} |\vc{D}| <  \infty$ and $\vc{\pi} |\vc{H}_{\bbA}| <  \infty$. Hence, it follows from \cite[Proposition~1.1]{Glyn96} that there exists some row vector $\vc{\beta}$ such that $\vc{D} = \vc{H}_{\bbA} + \vc{e}\vc{\beta}$. Since $\vc{\pi}\vc{D}=\vc{0}$, we have $\vc{\beta} = - \vc{\pi}\vc{H}_{\bbA}$, which leads to (\ref{eqn:D}). It also follows from  \cite[Proposition~1.1]{Glyn96} that two arbitrary solutions, denoted by $\vc{X}_1$ and $\vc{X}_2$, for the Poisson equation (\ref{Poisson-EQ-H}) with the constraint (\ref{Poisson-EQ-H-constraint}) satisfy  $\vc{X}_1 - \vc{X}_2 = \vc{e}\vc{\xi}$ for some row vector $\vc{\xi}$. Since $\vc{\pi}\vc{X}_1=\vc{\pi}\vc{X}_2 = \vc{0}$, 
\begin{align*}
\vc{0}
= \vc{\pi}(\vc{X}_1 - \vc{X}_2) = \vc{\pi} \vc{e}\vc{\xi} = \vc{\xi}.
\end{align*}
Therefore, $\vc{X}_1 = \vc{X}_2$, which implies that $\vc{D}$ is the unique solution for the Poisson equation (\ref{Poisson-EQ-H}) with the constraint (\ref{Poisson-EQ-H-constraint}). It has been proved that the statement (ii) holds if the statement (i) is true.

Next, we prove the statement (i). As its first step, we confirm that the statement (i) holds if the equivalence of (\ref{sum-pi(k)-u(k)-finite}) and (\ref{sum-pi(k,i)*E[T_{(0,j)}]-finite}) below is true:
\begin{align}
\sum_{k=0}^{\infty} \vc{\pi}(k) \vc{u}(k)
&< \infty,
\label{sum-pi(k)-u(k)-finite}
\\
\sum_{(k,i)\in\bbS} \pi(k,i)\EE_{(k,i)}[\tau_{(0,j)}] &< \infty
\quad \mbox{for some (and then every) $j \in \bbM_0$,}
\label{sum-pi(k,i)*E[T_{(0,j)}]-finite}
\end{align}
where $\tau_{(0,j)} = \inf\{n\in\bbN: (X_n,J_n) = (0,j)\}$ for $j \in \bbM_0$.
Assumption~\ref{assumpt-ergodic-MG1}, together with the aperiodicity of $\vc{P}$, ensures that the Markov chain $\{(X_n, J_n)\}$ is ergodic. Thus, the deviation matrix $\vc{D}$ exists if and only if (\ref{sum-pi(k,i)*E[T_{(0,j)}]-finite}) holds (see \cite[Theorem~4.1]{Cool02}). Additionally, Theorem~\ref{thm-u(k)} shows that (\ref{sum-pi(k)-u(k)-finite}) holds if and only if Assumption~\ref{assum-2nd-moment} does. Therefore, the equivalence of (\ref{sum-pi(k)-u(k)-finite}) and (\ref{sum-pi(k,i)*E[T_{(0,j)}]-finite}) ensures that Assumption~\ref{assum-2nd-moment} is equivalent to the existence of $\vc{D}$.

To complete the proof, we prove the equivalence of (\ref{sum-pi(k)-u(k)-finite}) and (\ref{sum-pi(k,i)*E[T_{(0,j)}]-finite}). For all $j \in \bbM_0$ and $(k,i) \in \bbS$, we have
\begin{align}
\EE_{(k,i)}[\tau_0]
&\le 
\EE_{(k,i)}[\tau_{(0,j)}]
\nonumber
\\
&= \EE_{(k,i)}[\tau_0] + \EE_{(k,i)}\!
\left[\one(J_{\tau_0} \neq j) \cdot \EE[ \tau_{(0,j)} - \tau_0 \mid J_{\tau_0} \neq j] 
\right]
\nonumber
\\
&= \EE_{(k,i)}[\tau_0] 
+ \sum_{s \in \bbM_0 \setminus \{j\}}
\PP(J_{\tau_0} = s \mid (X_0,J_0)=(k,i)) \cdot \EE_{(0,s)}[ \tau_{(0,j)}],
\label{eqn-191101-01}
\end{align}
where the last equality is due to the strong Markov property. Since $\{(X_n, J_n)\}$ is ergodic, there exists some constant $C>0$ such that
\[
\sup_{s,j \in \bbM_0}\EE_{(0,s)}[ \tau_{(0,j)}] < C.
\]
This bound and (\ref{eqn-191101-01}) yield
\begin{equation*}
\EE_{(k,i)}[\tau_0]
\le 
\EE_{(k,i)}[\tau_{(0,j)}]
\le \EE_{(k,i)}[\tau_0]  + C
\quad\mbox{for all $j \in \bbM_0$ and $(k,i) \in \bbS$},
\end{equation*}
which implies the equivalence of (\ref{sum-pi(k)-u(k)-finite}) and (\ref{sum-pi(k,i)*E[T_{(0,j)}]-finite}).

\section{Proof of Theorem~\ref{THM:bound-ol{pi}^{(N)}(k)}}\label{proof:THM:bound-ol{pi}^{(N)}(k)}

This appendix provides the proof of Theorem~\ref{THM:bound-ol{pi}^{(N)}(k)}, which requires the following lemma (its proof is given in Appendix~\ref{subsec:proof:LEM:pathiwse-ordering}).
\begin{LEM}\label{LEM:pathiwse-ordering}
For $k \in \bbZ_{[0,N-1]}$ and $j\in\bbM$,
\begin{align}
&
\sum_{\ell=k+1}^N
\EE_{(0,j)}^{(N)}\!\!\left[ 
\sum_{n=1}^{\tau_0^{(N)}} 
1_{(\ell,i)}(X_n^{(N)},J_n^{(N)})
\right]
\le  
\sum_{\ell=k+1}^{\infty}
\EE_{(0,j)}\!\!\left[ 
\sum_{n=1}^{\tau_0} 
1_{(\ell,i)}(X_n,J_n)
\right],
\label{eqn:210216-01}
\end{align}
where we use the simplified notation $\EE_{(k,i)}^{(N)}[\,\cdot\,] := \EE[\,\cdot\mid (X_0^{(N)}, J_0^{(N)})=(k,i)]$, in addition to the notation $\EE_{(k,i)}[\,\cdot\,]$ introduced in Section~\ref{subsec:mean_first passage_time}. 
\end{LEM}

In the following, we first prove Theorem~\ref{THM:bound-ol{pi}^{(N)}(k)} by using Lemma~\ref{LEM:pathiwse-ordering}, and then provide the proof of the lemma.

\subsection{Main body of the proof}

To prove (\ref{asymp-pi^{(N)}(k)}), it suffices to show that 
\begin{align}
\ol{\vc{\pi}}^{(N)}(k) 
\le (1 + o(N^{-1})) \ol{\vc{\pi}}(k),\qquad k \in \bbZ_{[0,N]},
\label{asymp-pi^{(N)}(k)-02}
\end{align}
where the term $1 + o(N^{-1})$ is independent of $k$. It follows from \cite[Theorem~10.0.1]{Meyn09} that, for $\ell \in \bbZ_{[0,N-1]}$ and $i \in \bbM_1$,
\begin{align}
\pi(\ell,i) &= 
\sum_{j \in \bbM_0}\pi(0,j)
\EE_{(0,j)}\!\!\left[ 
\sum_{n=1}^{\tau_0} 
1_{(\ell,i)}(X_n,J_n)
\right],
\label{eqn:pi(l,i)}
\\
\pi^{(N)}(\ell,i)
&= \sum_{j \in \bbM_0} \pi^{(N)}(0,j)
\EE_{(0,j)}^{(N)}\!\!\left[ 
\sum_{n=1}^{\tau_0^{(N)}} 
1_{(\ell,i)}(X_n^{(N)},J_n^{(N)})
\right].
\label{eqn:pi^{(N)}(l,i)}
\end{align}
It also follows from (\ref{eqn:210216-01}) and (\ref{eqn:pi^{(N)}(l,i)}) that, for $k \in \bbZ_{[0,N]}$ and $i \in \bbM_1$,
\begin{align}
\ol{\pi}^{(N)}(k,i)
&=
\sum_{j \in \bbM_0}
\pi^{(N)}(0,j)
\sum_{\ell=k+1}^N
\EE_{(0,j)}^{(N)}\!\!\left[ 
\sum_{n=1}^{\tau_0^{(N)}} 
1_{(\ell,i)}(X_n^{(N)},J_n^{(N)})
\right]
\nonumber
\\
&\le
\sum_{j \in \bbM_0}
\pi^{(N)}(0,j)
\sum_{\ell=k+1}^{\infty}
\EE_{(0,j)}\!\!\left[ 
\sum_{n=1}^{\tau_0} 
1_{(\ell,i)}(X_n,J_n)
\right].
\label{eqn:bound_pi^{(N)}(k,i)}
\end{align}
In addition, Theorem~\ref{thm-asymp-pi^(N)(k)} implies that $\vc{\pi}^{(N)}(0) \le (1 + o(N^{-1}))\vc{\pi}(0)$. Applying this to (\ref{eqn:bound_pi^{(N)}(k,i)}), and using (\ref{eqn:pi(l,i)}), we obtain
\begin{align*}
\ol{\pi}^{(N)}(k,i)
&\le
 (1 + o(N^{-1}))
\sum_{j \in \bbM_0}\pi(0,j)
 \sum_{\ell=k+1}^{\infty}
\EE_{(0,j)}\!\!\left[ 
\sum_{n=1}^{\tau_0} 
1_{(\ell,i)}(X_n,J_n)
\right]
\nonumber
\\
&=
 (1 + o(N^{-1}))
 \sum_{\ell=k+1}^{\infty}
\sum_{j \in \bbM_0}\pi(0,j)
\EE_{(0,j)}\!\!\left[ 
\sum_{n=1}^{\tau_0} 
1_{(\ell,i)}(X_n,J_n)
\right]
\nonumber
\\
&= (1 + o(N^{-1}) \sum_{\ell=k+1}^{\infty} \pi(\ell,i)
\nonumber
\\
&=
 (1 + o(N^{-1})) \ol{\pi}(k,i),\qquad k \in \bbZ_{[0,N]},~i\in\bbM_1,
\end{align*}
where the term $(1 + o(N^{-1}))$ is independent of $(k,i)$. Consequently, (\ref{asymp-pi^{(N)}(k)-02}) holds.

\subsection{Proof of Lemma~\ref{LEM:pathiwse-ordering}}\label{subsec:proof:LEM:pathiwse-ordering}

To prove Lemma~\ref{LEM:pathiwse-ordering}, we introduce a new stochastic process, and then prove a proposition (Proposition~\ref{PROP:Appendix-01}) on the new stochastic process. The proposition directly leads to the proof of Lemma~\ref{LEM:pathiwse-ordering}.

We construct a stochastic process $\{(\bv{X}_n^{(N)},\bv{J}^{(N)}_n);n\in\bbZ_+\}$ from the infinite-level M/G/1-type Markov chain $\{(X_n,J_n);n\in\bbZ_+\}$. Let $\{(\bv{X}_n^{(N)},\bv{J}_n^{(N)});n\in\bbZ_+\}$ denote a stochastic process such that
\begin{subequations}\label{defn:(ol{X}_n,ol{J}_n)}
\begin{alignat}{2}
\bv{X}_n^{(N)} &= 
\left\{
\begin{array}{lll}
X_0 \vmin N, & n=0,
\\
(\bv{X}_{n-1}^{(N)} + X_n - X_{n-1}) \vmin N, &
 n \in \bbN, &
\bv{X}_{n-1}^{(N)} \in \bbZ_{[0,N-1]},
\\
\bv{X}_{n-1}^{(N)} + (X_n - X_{n-1}) \vmin 0, &
 n \in \bbN, &
\bv{X}_{n-1}^{(N)} =N,
\end{array}
\right.
\label{defn:(ol{X}_n,ol{J}_n)-b}
\\
\bv{J}_n^{(N)} &= J_n, \qquad n \in \bbZ_+.
\label{eqn:bv{J}_n^{(N)}=bv{J}_n}
\end{alignat}
\end{subequations}
For later convenience, we refer to the value of $\{\bv{X}_n^{(N)}\}$ as {\it level}.

The following proposition implies that Lemma~\ref{LEM:pathiwse-ordering} holds.
\begin{PROP}\label{PROP:Appendix-01}
Let $\bv{\tau}_{0}^{(N)} = \inf\{n \in \bbN: \bv{X}_n^{(N)} = 0\}$, then the following are true.
\hfill
\begin{enumerate}
\item 
$\bv{X}_n^{(N)} \le X_n$ for all $n=0,1,\dots,\bv{\tau}_0^{(N)}$, and
 thus $\bv{\tau}_{0}^{(N)} \le \tau_0$.
\item The stopped process $\{(\bv{X}_n^{(N)},\bv{J}_n^{(N)}); 0 \le n \le \bv{\tau}_{0}^{(N)}\}$ is driven by the finite-level M/G/1-type transition probability matrix $\vc{P}^{(N)}$ until it reaches level zero, that is, while $ 0 \le n < \bv{\tau}_{0}^{(N)}$. 
\item
If $(\bv{X}_0^{(N)},\bv{J}_0^{(N)}) = (X_0^{(N)},J_0^{(N)})$, then $\bv{\tau}_{0}^{(N)}=\tau_{0}^{(N)}$ in distribution.
\end{enumerate}
\end{PROP}

\begin{proof}
The statement (i) is obvious, and therefore we prove the remaining two statements (ii) and (iii). By definition, the level increment of the stopped process $\{(\bv{X}_n^{(N)},\bv{J}_n^{(N)}); 0 \le n \le \bv{\tau}_{0}^{(N)}\}$ is equal to that of $\{(X_n,J_n); 0 \le n \le \bv{\tau}_{0}^{(N)}\}$ when the former is away from the upper boundary level $N$. Furthermore, the part  beyond level $N$ of the level increment is discarded; and once $\{(\bv{X}_n^{(N)},\bv{J}_n^{(N)})\}$ reaches level $N$, it keeps staying level $N$ with the substochastic matrix $\ol{\vc{A}}(-1)=\sum_{\ell=0}^{\infty}\vc{A}(\ell)$ until the infinite-level chain $\{(X_n,J_n)\}$ decreases its level with the substochastic matrix $\vc{A}(-1)$. Therefore, the statement (ii) holds (we can find a similar argument in \cite[Theorem~8.3.1]{Lato99}). In addition, for the proof of the statement (iii), we note that $\bv{X}_n^{(N)} = 1$ for $n=\bv{\tau}_{0}^{(N)}-1$. It thus follows from (\ref{defn:(ol{X}_n,ol{J}_n)}) and $\vc{B}(-1)\vc{e} = \vc{A}(-1)\vc{e}$ that
\begin{align*}
\PP(\bv{X}_{n+1}^{(N)} = 0 \mid \bv{X}_n^{(N)}=1,\bv{J}_n^{(N)} =i)
&=
\PP(X_{n+1} =  X_n - 1 \mid X_n \ge 1,J_n =i)
\nonumber
\\
&= (\vc{B}(-1)\vc{e})_i
\nonumber
\\
&= \PP(X_{n+1}^{(N)} = 0 \mid X_n^{(N)}=1,J_n^{(N)} =i).
\end{align*}
Consequently, the statement (ii) implies the statement (iii). The proof is completed.
\end{proof}

\section{Proof of Theorem~\ref{thm-asym}}\label{proof:thm-asym}

Theorem~\ref{thm-asym} is an immediate consequence of applying Lemmas~\ref{lem:main-01} and \ref{lem:main-02} to (\ref{eq-xN(k)-x(k)}). The two lemmas are given in Sections~\ref{proof:lem:main-01} and \ref{proof:lem:main-02}, respectively.
\begin{LEM}\label{lem:main-01}
If Assumptions~\ref{assumpt-ergodic-MG1}, \ref{assum-2nd-moment}, and \ref{assumpt-tail} hold, then
\begin{subequations}
\begin{alignat}{2}
\lim_{N\to\infty}
{1 \over \ol{F}(N)}
\left[
\vc{\pi}^{(N)}(0)
\sum_{n=N+1}^{\infty}\vc{B}(n)\vc{S}^{(N)}(n;k)
\right]
&=\vc{0}, & \qquad k &\in \bbZ_+,
\label{eqn:210216-02}
\\
\lim_{N\to\infty}
{1
\over 
\ol{F}(N)
}
\left[
\sum_{\ell=1}^N \vc{\pi}^{(N)}(\ell) 
\sum_{n=N+1}^{\infty} \vc{A}(n-\ell) \vc{S}^{(N)}(n;k)
\right]
&=\vc{0},& \qquad k &\in \bbZ_+,
\label{eqn:210216-03}
\end{alignat}
\end{subequations}
where $\vc{S}^{(N)}(n;k)$ is given in (\ref{def-S(N;n,k)}).
\end{LEM}

\begin{LEM}\label{lem:main-02}
If Assumptions~\ref{assumpt-ergodic-MG1}, \ref{assum-2nd-moment}, and \ref{assumpt-tail} hold, then
\begin{subequations}
\begin{align}
\lim_{N\to\infty}
{
\vc{\pi}^{(N)}(0)
\ool{\vc{B}}(N-1)\vc{e}
\over 
\ol{F}(N)
}
&= \vc{\pi}(0)\vc{c}_{B},
\label{eqn:210216-10}
\\
\lim_{N\to\infty}
\sum_{\ell=1}^N 
{
\vc{\pi}^{(N)}(\ell) 
\ool{\vc{A}}(N-\ell-1)\vc{e}
\over 
\ol{F}(N)
}
& = \ol{\vc{\pi}}(0)\vc{c}_{A}.
\label{eqn:210216-11}
\end{align}
\end{subequations}
\end{LEM}

\subsection{Proof of Lemma~\ref{lem:main-01}}\label{proof:lem:main-01}

We confirm that (\ref{eqn:210216-02}) and (\ref{eqn:210216-03}) hold if
\begin{subequations}
\begin{align}
\lim_{N\to\infty}
{
\vc{\pi}^{(N)}(0)\ol{\vc{B}}(N)\vc{e}
\over 
\ol{F}(N)
}
&=0,
\label{eqn:210216-04}
\\
\lim_{N\to\infty}
{
\sum_{\ell=1}^N \vc{\pi}^{(N)}(\ell) \ol{\vc{A}}(N-\ell)\vc{e}
\over 
\ol{F}(N)
}
&=0.
\label{eqn:210216-05}
\end{align}
\end{subequations}
Since $\vc{G}$ is stochastic, $\vc{G}^N = \ol{\bcalO}(1)$. It thus follows from (\ref{def-S(N;n,k)}) that, for each $k \in \bbZ_+$, there exists some finite $C_k > 0$ such that $|\vc{S}^{(N)}(n;k) | \le C_k \vc{e}\vc{e}^{\top}$ for all $n \ge N+1$ and $N \ge \max(k,1)$, which yields
\begin{subequations}
\begin{alignat}{2}
\sum_{n=N+1}^{\infty}\vc{B}(n) |\vc{S}^{(N)}(n;k)|
&\le C_k \ol{\vc{B}}(N) \vc{e}\vc{e}^{\top}, 
&\quad k & \in \bbZ_+,~N \ge \max(k,1),
\label{eqn:210316-01a}
\\
\sum_{n=N+1}^{\infty} \vc{A}(n-\ell) |\vc{S}^{(N)}(n;k)|
&\le C_k\ol{\vc{A}}(N-\ell)\vc{e}\vc{e}^{\top},
&\quad k & \in \bbZ_+,~N \ge \max(k,1).
\label{eqn:210316-01b}
\end{alignat}
\end{subequations}
Combining (\ref{eqn:210216-04}) and (\ref{eqn:210316-01a}) leads to (\ref{eqn:210216-02}); and combining (\ref{eqn:210216-05}) and (\ref{eqn:210316-01b}) leads to (\ref{eqn:210216-03}).

We first prove (\ref{eqn:210216-04}). Since $F \in \calL$, $\lim_{N\to\infty}\ol{F}(N-1) / \ol{F}(N)=1$. Therefore, using (\ref{defn:ool{A}(n)-ool{B}(n)}), (\ref{asymp-ool{B}(k)}), and Theorem~\ref{thm-asymp-pi^(N)(k)}, we obtain
\begin{align}
\lim_{N\to\infty}
{
\vc{\pi}^{(N)}(0)\ol{\vc{B}}(N)\vc{e}
\over 
\ol{F}(N)
}
&=
\lim_{N\to\infty}
{
\vc{\pi}^{(N)}(0) 
\left[ \ool{\vc{B}}(N-1)\vc{e} - \ool{\vc{B}}(N)\vc{e} \right]
\over 
\ol{F}(N)
}
\nonumber
\\
&=
\lim_{N\to\infty}
{
\vc{\pi}^{(N)}(0) 
\ool{\vc{B}}(N-1)\vc{e}
\over 
\ol{F}(N-1)
}
{\ol{F}(N-1) \over \ol{F}(N)}
-
\lim_{N\to\infty}
{
\vc{\pi}^{(N)}(0) 
\ool{\vc{B}}(N)\vc{e}
\over 
\ol{F}(N)
}
\nonumber
\\
&= 
\vc{\pi}(0)\vc{c}_{B} -\vc{\pi}(0)\vc{c}_{B}
= 0,
\label{eqn:210316-02}
\end{align}
which shows that (\ref{eqn:210216-04}) holds.

Next, we prove (\ref{eqn:210216-05}). Let $\Pi \ge 1$ and $\ol{A}$ denote independent integer-valued random variables. We then have
\begin{align}
\PP(\Pi + \ol{A} \ge N) = \PP(\Pi \ge 1)\PP(\ol{A} \ge N) 
+ \sum_{\ell=0}^{N-1} \PP(\Pi \ge N-\ell)\PP(\ol{A} = \ell),
\quad N \in \bbN. 
\label{eqn:220528-01}
\end{align}
By analogy with this equation, the following hold:
\begin{align}
\sum_{k=N}^{\infty} \sum_{\ell=1}^k 
\vc{\pi}^{(N)}(\ell) \ol{\vc{A}}(k-\ell)\vc{e}
&=
\ol{\vc{\pi}}^{(N)}(0) \ool{\vc{A}}(N-1)\vc{e}
+ \sum_{\ell=0}^{N-1} \ol{\vc{\pi}}^{(N)}(N - \ell - 1) \ol{\vc{A}}(\ell)\vc{e}.\label{eqn:210216-06a}
\\
\sum_{k=N}^{\infty}
\sum_{\ell=1}^k \vc{\pi}(\ell) \ol{\vc{A}}(k-\ell)\vc{e}
&=
\ol{\vc{\pi}}(0) \ool{\vc{A}}(N-1)\vc{e}
+ \sum_{\ell=0}^{N-1} \ol{\vc{\pi}}(N - \ell - 1) \ol{\vc{A}}(\ell)\vc{e}.
\label{eqn:210306-01}
\end{align}
Applying Theorem~\ref{THM:bound-ol{pi}^{(N)}(k)} to (\ref{eqn:210216-06a}) and using (\ref{eqn:210306-01}), we obtain
\begin{align}
&
\sum_{k=N}^{\infty}
\sum_{\ell=1}^k \vc{\pi}^{(N)}(\ell) \ol{\vc{A}}(k-\ell)\vc{e}
\nonumber
\\
&\le (1 + o(N^{-1}))
\left[
\ol{\vc{\pi}}(0) \ool{\vc{A}}(N-1)\vc{e}
+ \sum_{\ell=0}^{N-1} \ol{\vc{\pi}}(N - \ell - 1) \ol{\vc{A}}(\ell)\vc{e}
\right]
\nonumber
\\
&=
(1 + o(N^{-1}))
\sum_{k=N}^{\infty}
\sum_{\ell=1}^k \vc{\pi}(\ell) \ol{\vc{A}}(k-\ell)\vc{e}.
\label{eqn:210216-06b}
\end{align}
Furthermore, it follows from (\ref{asymp-ool{A}(k)}), Proposition~\ref{prop-Masu16}, and \cite[Proposition~A.3]{Masu11} that
\begin{align}
\lim_{N\to\infty}
\sum_{k=N}^{\infty}
\sum_{\ell=1}^k 
{
\vc{\pi}(\ell) \ol{\vc{A}}(k-\ell)\vc{e}
\over
\ol{F}(N)
}
=  { 
\vc{\pi}(0) \vc{c}_{B} + \ol{\vc{\pi}}(0) \vc{c}_{A}
\over
-\sigma 
} \vc{\varpi}
\ool{\vc{A}}(-1)\vc{e}
+ \ol{\vc{\pi}}(0)\vc{c}_A.
\label{eqn:210216-08}
\end{align}
Combining (\ref{eqn:210216-08}) and (\ref{eqn:210216-06b}) yields
\begin{align}
\limsup_{N\to\infty}
\sum_{k=N}^{\infty}
\sum_{\ell=1}^k 
{
\vc{\pi}^{(N)}(\ell) \ol{\vc{A}}(k-\ell)\vc{e}
\over
\ol{F}(N)
} 
&< \infty,
\label{eqn:210216-13}
\end{align}
and thus
\begin{align*}
\lim_{N\to\infty}
\sum_{\ell=1}^N
{
\vc{\pi}^{(N)}(\ell) \ol{\vc{A}}(N-\ell)\vc{e}
\over
\ol{F}(N)
} = 0,
\end{align*}
which shows that (\ref{eqn:210216-05}) holds. The proof is completed.

\subsection{Proof of Lemma~\ref{lem:main-02}}\label{proof:lem:main-02}

We prove only the second limit (\ref{eqn:210216-11}) because the first one (\ref{eqn:210216-10}) is implied by (\ref{eqn:210316-02}) in the proof of (\ref{eqn:210216-04}). Note that a similar equation to (\ref{eqn:220528-01}) holds:
\begin{align*}
\PP(\Pi + \ol{A} \ge N) = \PP(\Pi \ge N+1)\PP(\ol{A} \ge 0) 
+ \sum_{\ell=1}^N \PP(\Pi = \ell)\PP(\ol{A} \ge N - \ell),
\quad N \in \bbN. 
\end{align*}
From this equation, we have
\begin{align}
\sum_{k=N}^{\infty} \sum_{\ell=1}^k  
\vc{\pi}^{(N)}(\ell) \ol{\vc{A}}(k-\ell)\vc{e}
&= \ol{\vc{\pi}}^{(N)}(N) \ool{\vc{A}}(-1)\vc{e}
+ \sum_{\ell=1}^N \vc{\pi}^{(N)}(\ell) \ool{\vc{A}}(N-\ell-1)\vc{e}
\nonumber
\\
&= \sum_{\ell=1}^N \vc{\pi}^{(N)}(\ell) \ool{\vc{A}}(N-\ell-1)\vc{e},
\label{eqn:220528-02}
\end{align}
where the second equality is due to $\ol{\vc{\pi}}^{(N)}(N) =\sum_{k=N+1}^{\infty}\vc{\pi}^{(N)}(k) = \vc{0}$. Combining (\ref{eqn:220528-02}) and (\ref{eqn:210216-13}) yields
\begin{align*}
\limsup_{N\to\infty}
\sum_{\ell=1}^N 
{
\vc{\pi}^{(N)}(\ell) 
\ool{\vc{A}}(N-\ell-1)\vc{e}
\over 
\ol{F}(N)
}
< \infty.
\end{align*}
Therefore, using (\ref{asymp-ool{A}(k)}), Theorem~\ref{thm-asymp-pi^(N)(k)}, and the dominated convergence theorem, we obtain
\begin{align*}
&
\lim_{N\to\infty}
\sum_{\ell=1}^N 
{
\vc{\pi}^{(N)}(\ell) 
\ool{\vc{A}}(N-\ell-1)\vc{e}
\over 
\ol{F}(N)
}
\nonumber
\\
&\quad= 
\sum_{\ell=1}^{\infty} 
\lim_{N\to\infty}
{
\vc{\pi}^{(N)}(\ell) 
\ool{\vc{A}}(N-\ell-1)\vc{e}
\over 
\ol{F}(N)
}
\nonumber
\\
&\quad= 
\sum_{\ell=1}^{\infty} 
\lim_{N\to\infty}
\vc{\pi}^{(N)}(\ell) 
{
\ool{\vc{A}}(N-\ell-1)\vc{e}
\over 
\ol{F}(N -\ell-1)
}
{
\ol{F}(N -\ell-1)
\over 
\ol{F}(N)
}
\nonumber
\\
&\quad= 
\sum_{\ell=1}^{\infty} 
\vc{\pi}(\ell) 
\vc{c}_{A}
= \ol{\vc{\pi}}(0) \vc{c}_{A},
\end{align*}
which shows that (\ref{eqn:210216-11}) holds. The proof is completed.

%%%%%%%%%%%%%%%%%%%%%%%%%% Acknowledgments %%%%%%%%%%%%%%%%%%%%%%%%%%

\section*{Acknowledgments}
The research of Hiroyuki Masuyama was supported in part by JSPS KAKENHI Grant Number JP21K11770.

%%%%%%%%%%%%%%%%%%%%%%%%%%%%%%%%%%%%%%%%%%%%%%%%%%%%%%%%%%%%%%%%%%%%%%
%%%							Springer
%%%							Elsevier
%
%\bibliographystyle{elsarticle-harv}
%%%							Plain
\bibliographystyle{plain} % plain, alpha, abbrv, unsrt
%%%%%%%%%%%%%%%%%%%%%%%%%%%%%%%%%%%%%%%%%%%%%%%%%%%%%%%%%%%%%%%%%%%%%%
%\bibliography{hm2022_0527}

%%%%%%%%%%%%%%%%%%%%%%%%%%%%%%%%%%%%%%%%%%%%%%%%%%%%%%%%%%%%%%%%%%%%%%

\end{document}